\documentclass{amsart}
\usepackage{amssymb,amsmath}

\newtheorem{theorem}{Theorem}
\newtheorem{definition}{Definition}

\newtheorem{proposition}{Proposition}
\newtheorem{corollary}{Corollary}

\newcommand{\UP}{\overline{P}}
\newcommand{\LP}{\underline{P}}
\newcommand{\LQ}{\underline{Q}}

\newcommand{\LE}{\underline{E}}

\newcommand{\LU}{\underline{U}}

\newcommand{\rset}{\mathbb{R}}
\newcommand{\natset}{\mathbb{N}}
\newcommand{\intset}{\mathbb{Z}}
\newcommand{\reals}{\rset}
\newcommand{\naturals}{\mbox{I\hspace{-.1cm}N}}

\newcommand{\nset}{\naturals}

\newcommand{\dset}{\mathcal{D}}

\newcommand{\bset}{\mathcal{B}}
\newcommand{\cset}{\mathcal{C}}
\newcommand{\eset}{\mathcal{E}}
\newcommand{\fset}{\mathcal{F}}
\newcommand{\kset}{\mathcal{K}}

\newcommand{\hset}{\mathcal{H}}
\newcommand{\lset}{\mathcal{L}}

\newcommand{\pset}{\mathcal{P}}

\newcommand{\xset}{\mathcal{X}}

\newcommand{\setneg}[1]{#1^{c}}

\newcommand{\hqed}{\hfill\qed}

\newcommand{\xb}{\mathcal{X}_{B}}
\newcommand{\xomega}{\mathcal{X}_{\Omega}}

\newcommand{\fineproof}{\hfill \hqed}

\begin{document}



\title{Williams Coherence and Beyond}


\author{RENATO PELESSONI}

\address{DEAMS ``B. de Finetti''\\
University of Trieste\\
Piazzale Europa~1\\
I-34127 Trieste\\
Italy}
\email{renato.pelessoni@econ.units.it}

\author{PAOLO VICIG}

\address{DEAMS ``B. de Finetti''\\
University of Trieste\\
Piazzale Europa~1\\
I-34127 Trieste\\
Italy}
\email{paolo.vicig@econ.units.it}

\begin{abstract}
In this paper we discuss the consistency concept of Williams
coherence for imprecise conditional previsions, presenting a variant
of this notion, which we call W-coherence. It is shown that
W-coherence ensures important consistency properties and is quite
general and well-grounded. This is done comparing it with
alternative or anyway similar known and less known consistency
definitions. The common root of these concepts is that they
variously extend to imprecision the subjective probability approach
championed by de Finetti. The analysis in the paper is also helpful
in better clarifying several little investigated aspects of these
notions.
\end{abstract}

\maketitle

\section*{Acknowledgement}
*NOTICE: This is the authors' version of a work that was accepted for publication in the International Journal of Approximate Reasoning. Changes resulting from the publishing process, such as peer review, editing, corrections, structural formatting, and other quality control mechanisms may not be reflected in this document. Changes may have been made to this work since it was submitted for publication. A definitive version was subsequently published in the International Journal of Approximate Reasoning, vol. 50, issue 4, April 2009, doi:10.1016/j.ijar.2008.10.002

$\copyright$ Copyright Elsevier

http://www.sciencedirect.com/science/article/pii/S0888613X08001643

\section{Introduction}
\label{sec:introduction}
Quite recently, P.M. Williams' 1975 seminal paper \emph{Notes on conditional previsions} was published in a slightly revised
version \cite{wil07}, preceded by an introductory paper discussing basic aspects and historical motivations for
his work \cite{vic07}.
This fact confirms that Williams' ideas on coherence still play a very important role in the theory of conditional imprecise previsions.
In the past, they influenced the more widespread theory developed by Walley \cite{wa91}.
Williams coherence was also directly used in some papers to achieve results in different areas,
including epistemic independence \cite{vic00},
problems of checking consistency for conditional imprecise probabilities \cite{wa04},
consistency for unbounded random variables \cite{tro06}.
Yet, certainly also because of its overall limited diffusion in the scientific community,
several aspects of Williams coherence are still little explored.

A basic motivation for studying Williams coherence is its
generality: in the version we present in the paper,
it extends to a very broad conditional setting Walley's
(unconditional) coherence, which already encompasses as special cases
several uncertainty theories ($2$-monotone probabilities, precise probabilities, belief
functions, possibility/necessity measures, coherent risk
measures,...) applied in many different areas, from artificial
intelligence to statistics or risk measurement. Thus extensions of such theories to
conditional frameworks can be accommodated into
Williams coherence, exploiting hence the results already established
for it. In many cases, these problems have been so far little
investigated; for instance, much work remains to be done in the area of
measuring conditional risks. Williams coherence is not the only way
of extending Walley's (unconditional) coherence, but it is a very general
and (perhaps) immediate one; it is anyway important to weigh
pros and cons in choosing which coherence notion should be used.
This evaluation affects various issues, some more familiar (like the
validity of the envelope theorem), other ones generally less familiar (like
the problem of non-conglomerability).

The main purpose of this paper (extending earlier results in \cite{pel07}, Section 3)
is to investigate more closely the role of Williams coherence,
comparing it with the nearest consistency concepts that have been developed in the literature.
Since Williams' work was inspired by de Finetti's ideas,
these concepts are among those following and generalising the subjective probability approach to uncertainty.
We supply some historical information on this in Section \ref{sec:historical_note}.
The hints presented there are historically not exhaustive,
being limited to key contributions and under the perspective of studying Williams coherence,
but they let us mention a few important properties regarding all of these concepts,
which form a basis for making comparisons among them in the paper.
Section \ref{sec:other_preliminaries} contains other preliminary definitions and notions.

We investigate Williams coherence with a progressively larger perspective throughout the paper.
We start in Section \ref{sec:about_williams_definition} by discussing a nimbler variant for it,
called W-coherence (already defined in Section \ref{sec:other_preliminaries}), which is adopted in the sequel.
It generalises Walley coherence for unconditional lower previsions.
In Section \ref{sec:unconditional_conditional}, we discuss potential generalisations of other
unconditional coherence concepts,
focusing in particular on a little explored definition \cite{art99,maa03}.
We supply an interesting interpretation for it, showing that it has no straightforward extension to conditionals,
but deriving some conditions that are either necessary or sufficient for W-coherence.
In Section \ref{sec:alternative_concepts_coherence} W-coherence is compared with alternative views of
conditional coherence developed by Walley \cite{wa91},
proving in particular its equivalence with separate coherence
(when they are comparable, since separate coherence is less general).
The comparison is continued in Section \ref{sec:non_conglomerability},
discussing non-conglomerability and the different treatment of this property in Williams' and Walley's approaches.
Concepts related to W-coherence are discussed in Section \ref{sec:beyond_Williams_coherence}.
In particular, in Section \ref{sec:Williams_de_Finetti} we explore a notion intermediate between W-coherence
and dF-coherence, showing its little significance,
while in Section \ref{sec:aul_condition} we discuss which condition of avoiding loss - type should be
appropriate when adopting W-coherence.
It is shown that using a certain concept of avoiding uniform loss some seemingly inconsistent features of W-coherence
pointed out by Walley can be justified.
In Section \ref{sec:centered_convexity} we discuss centered convexity,
a relevant concept, (moderately) weaker than W-coherence.
We point out that centered convexity shares desirable properties with W-coherence,
even though this is true at a lesser extent as far as envelope theorems are concerned.
Our conclusions on the role of W-coherence in imprecise probability theory are contained
in Section \ref{sec:conclusions}.

\section{Preliminary Issues}
\label{sec:preliminaries}
We recall first a few notions concerning the description of uncertainty.
Following \cite{cri96,def74} and others, we use the logical notation to operate with events.
This originates from observing that events are described by propositions of classical logic,
and actually a formal definition of events and conditional events in these terms was given in \cite{cri96,cri06}.
We write $B$ for both an event $B$ and its indicator function $|B|$ (de Finetti's convention),
appearing from the context which of the two meanings is intended.

A \emph{bounded} random variable\footnote{Also called gamble in \cite{wa91} or bounded random quantity in \cite{wil07},
whilst random quantities can be unbounded in \cite{def74}.}
$X$ is \emph{represented} by a map $X:\bset\rightarrow\rset$, where $\bset$ is a partition of (non-impossible) events.
A possible value of $X$, $X(\omega)$, corresponds to each $\omega\in\bset$,
which does not mean that $\bset$ is unique.
For instance, the partition whose generic event is `$X=x$' will do
if we describe $X$ alone,
but a more refined partition is needed to describe two or more random variables simultaneously.
In classical probability theory, a unique \emph{fixed} partition (called $\Omega$ there, while we reserve the symbol $\Omega$
for the sure event), large enough to describe what matters, is employed.
This is not necessary in general (cf. the discussion in \cite{wa91}, Section 2.1.4) and will not be pursued here.

When conditioning on some non-impossible event $B$, the conditional random variable $X|B$ may be represented by
$X_{B}:\bset|B\rightarrow\rset$,
where the elements of the conditional partition $\bset|B$ are obtained replacing each $\omega\in\bset$
with the conditional event $\omega|B$, and discarding those $\omega|B$ which turn out to be impossible,
conditional on $B$ (i.e., such that assuming $B$ true implies that $\omega$ is false).
After this is done, $X_{B}(\omega|B)=X(\omega)$ holds.
In the special case that $\bset=\{\Omega\}$, we reobtain unconditional random variables ($X|\Omega=X$).

The supremum $\sup(X|B)$ of $X|B$ may be computed as
$\sup_{\omega\Rightarrow B}X(\omega)$ (in the set-theoretic language: $\sup_{\omega\in B}X(\omega)$).

When working with conditional random variables, we shall sometimes employ the equality
\begin{eqnarray}
\label{eq:equivalence_conditional_function}
f(X_{1},\ldots,X_{n})|B=f(X_{1}|B,\ldots,X_{n}|B)
\end{eqnarray}
where $f$ is any real function,
returning the random variable $f(X_{1},\ldots,X_{n})$ as a function of $X_{1},\ldots,X_{n}$ \cite{cri06}.
A typical case we will consider in the paper is $f=G$, where $G$ is a `gain'.

In the rest of the paper, the domain of the uncertainty measures considered is usually termed $\dset$.
Precisely, $\dset$ is an \emph{arbitrary} (non-empty) set of \emph{bounded} random variables,
or more generally of bounded conditional random variables.
$\dset$ may contain conditional events too,
corresponding to those $X|B\in\dset$ such that $X$ is the indicator of some event,
or events when further $B=\Omega$.

A \emph{lower prevision} $\LP$ on $\dset$ is a map $\LP:\dset\rightarrow\rset$.
An upper prevision $\UP$ may be defined through the equality $\UP(-X|B)=-\LP(X|B)\ \forall X|B\in\dset$,
which always lets us refer to either lower or upper previsions only.
A \emph{precise} prevision $P$ corresponds to the special case $\UP(X)=\LP(X)=P(X)$.

\subsection{A Historical Note}
\label{sec:historical_note}
We shall deal in this paper with several notions of `coherence', or weaker concepts.
Their forerunner is de Finetti's coherence for (unconditional) precise previsions \cite{def74}:
\begin{definition}
\label{def:coherent_prevision}
$P:\dset\rightarrow\reals$ is a
\emph{dF-coherent} precise \emph{prevision} on $\dset$
iff, $\forall n,m \in\natset$, $\forall\ X_1,\ldots,X_n,Y_{1},\ldots,Y_{m} \in
\dset$, $\forall\ s_1,\ldots,s_n\geq 0$, $\forall\ r_1,\ldots,r_m\geq 0$,
defining $G=\sum_{i=1}^{n}s_{i}(X_i-P(X_i))-\sum_{j=1}^{m}r_{j}(Y_j-P(Y_j))$,
it holds that $\sup G\geq 0$.
\end{definition}
This definition includes that of \emph{dF-coherent (precise) probability}
as a special case,
when all the random variables in $\dset$ are (indicators of) events.
If further $\dset$ is an algebra, dF-coherent probabilities coincide with finitely additive probabilities.

The notion of dF-coherent precise prevision is also closely related to that of \emph{expectation}:
whenever an expectation $E(X)$ is assessed for $X$,
then $E(X)$ is also its only dF-coherent prevision.
However, whenever $P(X)$ is assessed $E(X)$ is not necessarily defined,
because no probability on all events `$X\leq x$' must be preliminarily elicited in order to define $P(X)$.

Although de Finetti did not develop extensively a theory of conditional previsions,
nor was he much concerned with imprecise previsions, several features in his approach were influential
also in most later generalisations.
We mention the following basic facts, referring to a generic, not specified `consistency' property of
(precise or imprecise) previsions.
\begin{itemize}
\item[A)] Previsions are announced on an \emph{arbitrary} (non-empty) set of random variables $\dset$,
and consequently the definition of their consistency is structure-free.
\item[B)] An \emph{extension theorem} ensures that a consistent prevision can be extended on \emph{any}
$\dset^{\prime}\supset\dset$, so that the extension preserves the same type of consistency on $\dset^{\prime}$.
\item[C)] Consistent previsions have a behavioural interpretation in some idealised betting scheme.
\end{itemize}
When `consistency' is replaced by `dF-coherence', A), B) and C) are satisfied.\footnote{The dF-coherent extension is generally not unique.
In the special case that the events `$X\leq x$' are in $\dset$ $\forall x\in\rset$,
while $X\notin\dset$,
the dF-coherent extension on $\dset\cup\{X\}$ is unique,
and as mentioned above coincides with $E(X)$.

}
Concerning the betting scheme, the random variable $G$ in Definition \ref{def:coherent_prevision}
is the \emph{gain} from a bet made up of $n+m$ elementary bets,
$n$ `in favour of'  $X_{1},\ldots,X_{n}$ (the bettor is willing to pay
$s_{i}P(X_{i})$ for receiving $s_{i}X_{i}$, $i=1,\ldots,n$),
$m$ `against' $Y_{1},\ldots,Y_{m}$ (the bettor receives $r_{j}P(Y_{j})$ to sell $r_{j}Y_{j}$, $j=1,\ldots,m$).
The definition of dF-coherence requires that, whatever is the bet, the gain cannot be negative and bounded away from $0$.
DF-coherent previsions are linear and homogeneous, if the relevant quantities are in the domain $\dset$:
\begin{eqnarray}
\label{eq:linearity}
P(\alpha X+\beta Y)=\alpha P(X)+\beta P(Y).
\end{eqnarray}
When adding the constraint $m\leq 1$ in Definition \ref{def:coherent_prevision}, we obtain Walley's definition of
\emph{coherence} for lower previsions:
\begin{definition}
\label{def:coherent_lower_prevision}
$\LP:\dset\rightarrow\reals$ is a
\emph{coherent lower prevision} on $\dset$
iff, for all $n\in\natset$, $\forall\ X_0,X_1,\ldots,X_n \in
\dset$, $\forall\ s_0,s_1,\ldots,s_n\geq 0$,
defining $G=\sum_{i=1}^{n}s_{i}(X_i-\LP(X_i))-s_0(X_0-\LP(X_0))$,
it holds that $\sup G\geq 0$.
\end{definition}
Again, items A), B) and C) above are satisfied by this definition.
It has a well-known behavioural interpretation,
discussed in \cite{wa91}.
Some consequences of this interpretation,
not all highlighted in \cite{wa91},
may better stress the behavioural difference with Definition \ref{def:coherent_prevision}.
Precisely, $\LP(X)$ is an agent's supremum buying price for $X$, and $G$ is
the agent's \emph{gain} resulting from her/his buying $s_{i}X_{i}$,
for $i=1,\ldots,n$, and selling $s_{0}X_{0}$.
Coherence implies,
writing the last term in $G$ as $s_{0}\LP(X_{0})-s_{0}X_{0}$,
that the agent may be forced to accept \emph{at most one} of her/his supremum buying prices,
$s_{0}\LP(X_{0})$, as an infimum selling price for $s_{0}X_{0}$.
The restriction `at most one' does not apply to Definition \ref{def:coherent_prevision},
because $m$ there may be any natural number.
In general, we shall say that the agent \emph{bets} on (in favour or against)
$X$ with \emph{stake} $s$.

With imprecise previsions, there is a fourth property that we will consider:
\begin{itemize}
\item[D)] consistent imprecise previsions \emph{are characterised} by some envelope theorem.
\end{itemize}
Generally speaking, envelope theorems relate a function in a certain set $\fset$ to a set $\pset$ of other functions
with well specified features.
These theorems either ensure that by performing the (pointwise) infimum or supremum on the elements of $\pset$
we get a function $f\in \fset$, or else guarantee that every $f\in \fset$ may be expressed as an infimum or
supremum over some set $\pset$, or both (thus characterising the functions in $\fset$).
Envelope theorems may be found in many different research areas, like for instance cooperative games
\cite{sch72}.
They are important because:
\begin{itemize}
\item[$\bullet$] they ensure an often simple way of assigning a function $f$ with the desired consistency properties;
\item[$\bullet$] when being also characterisation theorems,
they allow an alternative, indirect definition and interpretation of the functions in $\fset$ by means of sets of the
(usually simpler) functions in $\pset$.
Moreover, they allow proving properties of the functions in $\fset$ using known results about the functions in $\pset$.
\end{itemize}
Coherent lower previsions ensure property D): a real function $\LP$ is a coherent lower prevision over $\dset$
if and only if $\LP(X)=\inf_{P\in\pset}\{P(X)\}, \forall X\in\dset$ ($\inf$ is attained),
where $\pset$ is a set of dF-coherent precise previsions $P$ dominating $\LP$ on $\dset$,
i.e. $P(X)\geq\LP(X)\ \forall X\in\dset, \forall P\in\pset$ \cite{wa91}.

Various generalisations of dF-coherence to conditional (precise or imprecise) previsions have been proposed.
The adherence of some of them to A), B) and D) will be discussed throughout the paper.
As for C), all of them have some behavioural interpretation. This aspect will therefore be just outlined.
In particular, dF-coherence for conditional (precise) previsions was developed in the eighties in \cite{hol85,reg85},
obeying the requirements A), B), C) above.
\begin{definition}
\label{def:conditional_precise_prevision}
$P:\dset\rightarrow\rset$ is a \emph{dF-coherent conditional (precise) prevision} on $\dset$ iff,
for all $n, m\in\natset$, $\forall\ X_1|B_1,\ldots,X_n|B_n,Y_1|C_1,\ldots,Y_m|C_m \in \dset$,
$\forall\ s_i\geq 0\ (i=1,\ldots,n)$,  $\forall\ r_j\geq 0\ (j=1,\ldots,m)$,
defining $G=\sum_{i=1}^{n}s_{i}B_{i}(X_i-P(X_i|B_i))-\sum_{j=1}^{m}r_{j}C_{j}(Y_j-P(Y_j|C_j))$,
$B=\bigvee_{i=1}^{n}B_i\vee\bigvee_{j=1}^{m}C_j$,
it holds that $\sup(G|B)\geq 0$.
\end{definition}
Here the gain is $G|B$, a conditional random variable itself.
Conditioning on $B$ has the meaning of considering only those values for $G$
when at least one of $B_{1},\ldots,B_{n},C_{1},\ldots,C_{m}$ is true.
Property (\ref{eq:linearity}) generalises to
\begin{eqnarray}
\label{eqn:linearity_coherence}
P(\alpha X+\beta Y|B)=\alpha P(X|B)+\beta P(Y|B).
\end{eqnarray}

Coherence concepts for conditional imprecise previsions were given by Walley \cite{wa91},
see Section \ref{sec:alternative_concepts_coherence}.
But the earliest proposal was that of Williams \cite{wil07} in 1975.
His work had a limited diffusion in those years, but influenced Walley's work and
contained \emph{in nuce} several fundamental results in the theory of imprecise probabilities \cite{vic07}.

\subsection{W-coherence and Other Preliminaries}
\label{sec:other_preliminaries}
In a conditional environment, we adopt the following generalisation of Definition \ref{def:coherent_lower_prevision}
to define a \emph{coherent} lower prevision $\LP(\cdot|\cdot)$:
\begin{definition}
\label{def:conditional coherent lower prevision}
$\LP:\dset\rightarrow\reals$ is a \emph{coherent conditional lower prevision} on $\dset$ iff, for all $n\in\natset$,
$\forall X_0|B_0,\ldots,X_n|B_n\in\dset$, $\forall\ s_0,s_1,\ldots,s_n$ real and non-negative,
defining $B=\bigvee_{i=0}^{n} B_{i}$ and
$G=\sum_{i=1}^{n}s_{i}B_{i}(X_{i}-\LP(X_{i}|B_{i}))-s_{0}B_{0}(X_{0}-\LP(X_{0}|B_{0}))$,
$\sup(G|B)\geq 0$.
\end{definition}
It is easy to realise that we would get an equivalent definition (adopted in \cite{wa04}) by replacing
$G|B$ with $G|S$, where the \emph{support} $S$ is defined as $S=\bigvee\{B_{i}:s_{i}\neq 0, i=0,\ldots,n\}$.

Throughout the paper, Definition \ref{def:conditional coherent lower prevision} will be referred to
as Williams coherence, or \emph{W-coherence} or simply coherence when unambiguous,
but as we will explain in Section \ref{sec:about_williams_definition}, it is actually a structure-free version
of the original Williams coherence.

A weaker notion than W-coherence is that of lower prevision that \emph{avoids uniform loss} \cite{wa04},
recalled in Section \ref{sec:aul_condition}.
In the unconditional environment it is termed condition of \emph{avoiding sure loss} and is defined in \cite{wa91},
Section 2.4.4 (a).

A further consistency notion, \emph{centered convexity} \cite{pel03bis, pel05,pel05bis}, is weaker than coherence,
but sufficiently stronger than the conditions of avoiding sure or uniform loss to allow for
interesting properties and applications (for instance, in risk measurement \cite{pel05}).
Its relationship with W-coherence is discussed in Section~\ref{sec:centered_convexity}.

Formally, the definition of \emph{convex lower prevision} is obtained from Definition \ref{def:coherent_lower_prevision}
and Definition \ref{def:conditional coherent lower prevision} by introducing just the extra \emph{convexity constraint}
$\sum_{i=1}^{n}s_{i}=s_{0}\ (>0)$ and eventually by further imposing (this is not restrictive) that $s_{0}=1$
\cite{pel03bis, pel05}.
Again, we could equivalently condition $G$ on its support $S$ rather than on $B$, as done in \cite{pel05,pel05bis}.
\emph{Centered} convexity requires in addition that ($0\in\dset$ and) $\LP(0)=0$ in the unconditional case,
and further that $\forall X|B\in\dset$, $0|B\in\dset$ and $\LP(0|B)=0$ in the conditional case
(cf. Definition \ref{def:conditional convex lower prevision}).
Centering is quite a natural requirement: non-centered convex previsions have rather weak consistency properties
(see also Footnote \ref{foo:conceptual_difference}),
but special instances of them may be found in the risk literature (cf. \cite{pel05}).

Let $\LP$ be a lower prevision defined on an arbitrary set $\dset$.
Following B) of Section \ref{sec:historical_note},
any consistency condition satisfied by $\LP$ should guarantee
that there exists an extension of $\LP$ satisfying the same
condition on any $\dset^{\prime}\supset\dset$.
If such an extension is not unique, its vaguest or least-committal one, if existing, has a special importance.
This peculiar extension is the \emph{natural extension} $\LE$ in the case of coherent or,
when conditioning, W-coherent previsions \cite{vic07,wa91}, the \emph{convex} natural extension $\LE_{c}$
for centered convex (unconditional or conditional) previsions \cite{pel03bis, pel05}.
The natural or convex natural extensions always exist for these consistency notions, not necessarily with other
ones, like Walley coherence in \cite{wa91}, Section 7.1.4 (b), or non-centered convexity.

\section {Coherence Concepts of Williams and Others}
\label{sec:Williams_others_coherence}
\subsection{About Williams' Definition}
\label{sec:about_williams_definition}
Williams' original definition (\cite{wil07}, Definition 1) differs formally from our definition of
W-coherence. One reason is that it refers to upper rather than lower previsions, but this is unimportant, since using
the conjugacy relation $\UP(-X|B)=-\LP(X|B)$ our condition $\sup (G|B)\geq 0$ corresponds
exactly to his inequality in $(A^{*})$ of \cite{wil07}.
The true difference is that his notion is not completely structure-free,
as it asks in particular that, for every conditioning event $B$,
the set $\xb=\{X: X|B \in \dset\}$ is a linear space.
It follows for instance that Williams' definition does not formally generalise Walley coherence for unconditional
previsions (our Definition \ref{def:coherent_lower_prevision}), which is structure-free:
when $B=\Omega$ for all $X|B\in\dset$, the set of all $X$ is constrained to form a linear space $\xomega$.
On the contrary, Definition \ref{def:conditional coherent lower prevision} is in particular a generalisation of
Walley's unconditional coherence and appears to be, in general, nimbler.
The fundamental link between the two versions of Williams coherence is ensured by the following \emph{extension theorem}.
\begin{proposition}
\label{pro:extension_W_coherent}
If $\LP:\dset\rightarrow\rset$ is W-coherent on $\dset$ (according to Definition \ref{def:conditional coherent lower prevision}),
it has a W-coherent extension on \emph{any} $\dset^{\prime}\supset\dset$.
\end{proposition}
Although we are not aware of any published proof for this
proposition, nevertheless it should be regarded as essentially
known. In fact, it can be proven by adapting the proofs concerning
the convex natural extension in \cite{pel05}, thus proving that
there \emph{always} exists the natural extension of a W-coherent
lower prevision on any $\dset^{\prime}\supset\dset$.
A proof of this kind is given in the Appendix, for the sake of completeness.
Alternatively,
the historically older scheme of de Finetti's extension theorem can be followed, with
suitable (but basically minor) modifications. After de Finetti's
path-breaking proof concerning precise (unconditional) previsions in
\cite{def49}, this scheme was employed in several generalisations
(see e.g. \cite{arm89, cri06,hol85}).
In the version for W-coherence, its two-step proof shows in the
first step that there exist W-coherent extensions on
$\dset^{\prime}=\dset\cup\{X|B\}$, $\forall X|B$, while the second
step generalises the proof to any $\dset^{\prime}$ using Zorn's
lemma or equivalent results. A by-product of the first step is that
the set of admissible W-coherent extensions on $X|B$ is proved to be
a closed interval. Its lower endpoint is the \emph{natural
extension} $\LE(X|B)$, while the upper endpoint is the \emph{upper
extension} $\LU(X|B)$ of $\LP$.
Thus, the scheme of de Finetti's extension theorem does not emphasise the role of the natural extension,
but rather treats the natural and upper extension in a symmetric way.

As an important implication of Proposition \ref{pro:extension_W_coherent} in our framework,
when $\dset$ in Definition \ref{def:conditional coherent lower prevision} does not meet the structure requirements
in Williams' definition it is always possible to coherently extend $\LP$ on a set $\dset^{\prime}$ such that these requirements hold,
and there the two notions of coherence coincide.
It follows that W-coherent lower previsions have all the properties established for Williams coherence in \cite{wil07},
including the important \emph{envelope theorem}, stating that $\LP$ is coherent on $\dset$ if and only if
\begin{eqnarray*}
\LP(X|B)=\inf_{P\in\pset}P(X|B), \forall X|B\in\dset
\end{eqnarray*}
where $\pset$ is a set of dF-coherent precise previsions $P(\cdot|\cdot)$ dominating $\LP(\cdot|\cdot)$ on $\dset$
($\forall P\in\pset$, $P(X|B)\geq\LP(X|B), \forall X|B\in\dset$).
Note that $\inf$ is attained.

\subsection{From Unconditional to Conditional Coherence}
\label{sec:unconditional_conditional}
As we have already pointed out, Definition \ref{def:conditional coherent lower prevision} of W-coherence generalises
Walley coherence for unconditional previsions (Definition \ref{def:coherent_lower_prevision}).
But other known definitions are equivalent to Definition \ref{def:coherent_lower_prevision}.
An interesting issue is therefore: why not rather generalise them in a conditional environment?
An answer is that Definition \ref{def:coherent_lower_prevision} seems more appropriate
for further generalisations.

The matter is relatively simple and well known if we consider a version of coherence,
equivalent to Definition \ref{def:coherent_lower_prevision},
obtained by restricting
the stakes $s_{0},\ldots,s_{n}$ to be integers (this is Walley's Definition 2.5.1 in \cite{wa91}).
The constraint on the integer stakes can be adopted in a conditional environment too,
for W-coherence as well as for some other consistency notions we discuss in this paper, obtaining equivalent formulations.
However, considering integer combinations only is not enough when the random variables are unbounded,
even in the unconditional case, as shown in \cite{tro02}.
We are not dealing with unbounded random quantities here,
yet in view of (potentially) pursuing the utmost generality,
we prefer not to impose the integer stakes constraint.

The situation is more complex, and definitely less explored, when turning to the following less used definition,
which is known to be equivalent to Definition~\ref{def:coherent_lower_prevision}:
\begin{definition}
\label{def:equivalent_definition_coherence}
$\LP:\dset\rightarrow\reals$ is a coherent lower prevision on $\dset$ iff,
for all $n\in\natset$, $\forall\ X_0,X_1,\ldots,X_n \in
\dset$, $\forall\ r_1,\ldots,r_n\geq 0$,
$\forall\ \mu_{0}\in\rset$ such that $X_{0}\geq\sum_{i=1}^{n}r_{i}X_{i}+\mu_{0}$,
it holds that
$\LP(X_{0})\geq\sum_{i=1}^{n}r_{i}\LP(X_{i})+\mu_{0}$.
\end{definition}
Definition \ref{def:equivalent_definition_coherence}
has a curious story: not mentioned explicitly in \cite{wa91},
although following directly from results established there,
it appears in \cite{art99}, but without being related to coherence for imprecise previsions,
which was later done in \cite{maa03}.
To the best of our knowledge,
Definition \ref{def:equivalent_definition_coherence} has not been given a clear behavioural interpretation yet,
nor has its potential generalisation to a conditional environment been explored.
We tackle these issues in this section.

As a first step, we rewrite the condition in Definition \ref{def:equivalent_definition_coherence},
that is,
\begin{eqnarray}
\label{eq:rewrite_alternative_coherence}
X_{0}\geq\sum_{i=1}^{n}r_{i}X_{i}+\mu_{0}\Rightarrow\LP(X_{0})\geq\sum_{i=1}^{n}r_{i}\LP(X_{i})+\mu_{0},
\end{eqnarray}
in an equivalent form.
Multiply for this the inequalities in (\ref{eq:rewrite_alternative_coherence}) by $s_{0}>0$,
let $s_{i}=r_{i}s_{0}$ $(i=1,\ldots,n)$, $\lambda_{0}=\mu_{0}s_{0}$ and perform the infimum in the first
inequality to obtain:
\begin{eqnarray}
\label{eq:rewrite_alternative_coherence_bis}
\lambda_{0}\leq\inf(s_{0}X_{0}-\sum_{i=1}^{n}s_{i}X_{i})\Rightarrow\lambda_{0}\leq s_{0}\LP(X_{0})-\sum_{i=1}^{n}s_{i}\LP(X_{i}).
\end{eqnarray}
If we define $I=-s_{0}X_{0}+\sum_{i=1}^{n}s_{i}X_{i}$, $E=-s_{0}\LP(X_{0})+\sum_{i=1}^{n}s_{i}\LP(X_{i})$,
(\ref{eq:rewrite_alternative_coherence_bis}) is rewritten as
\begin{eqnarray}
\label{eq:rewrite_alternative_coherence_ter}
\lambda_{0}\leq\inf(-I)\Rightarrow\lambda_{0}\leq-E.
\end{eqnarray}
Let us now come to the behavioural interpretation of Definition \ref{def:equivalent_definition_coherence}.
For any given bet on $X_{0},\ldots,X_{n}$ with stakes $s_{0},\ldots,s_{n}$,
$I$ is the bettor's overall \emph{income} ensuing from the bet,
while $E$ is her/his \emph{expense} for betting.
Note that $I$ is random, while $E$ is not.
From (\ref{eq:rewrite_alternative_coherence_ter}),
Definition \ref{def:equivalent_definition_coherence} asks as a necessary and sufficient condition for coherence
that $\inf(-I)\leq-E$, i.e. that $\sup(I)\geq E$, for \emph{any} bet.
This is a reasonable requirement: it does not hold iff $\sup I<E$ for \emph{some} bet,
and this means that a specific bet can be arranged whose ensuing gain $G=I-E$ is strictly negative and bounded away from zero
whatever happens, and the bettor suffers from a sure loss.
It is clear then that Definitions \ref{def:coherent_lower_prevision} and \ref{def:equivalent_definition_coherence}
are equivalent:
they both require that no bet must be such that $\sup G<0$.

The above interpretation also suggests a way to explore extensions of Definition \ref{def:equivalent_definition_coherence}
in a conditional framework.
Rewrite for this the gain $G$ in Definition \ref{def:conditional coherent lower prevision} highlighting the expense and
income terms.
We have $I=-s_{0}B_{0}X_{0}+\sum_{i=1}^{n}s_{i}B_{i}X_{i}$,
$E=-s_{0}B_{0}\LP(X_{0}|B_{0})+\sum_{i=1}^{n}s_{i}B_{i}\LP(X_{i}|B_{i})$
and the condition $\sup(G|B)\geq 0$ in Definition \ref{def:conditional coherent lower prevision} is written as
\begin{eqnarray}
\label{eq:extension_definition_coherence}
\sup(I-E|B)\geq 0.
\end{eqnarray}
The following Proposition is fundamental for discussing the potential generalisations of Definition \ref{def:equivalent_definition_coherence}.
\begin{proposition}
\label{pro:generalisation_coherence_definition}
Consider, as in Definition \ref{def:conditional coherent lower prevision}, a bet on $X_{0}|B_{0},\ldots,X_{n}|B_{n}$ with stakes
$s_{0},\ldots,s_{n}$, respectively, and define $B=\bigvee_{i=0}^{n}B_{i}$.
Let $\lambda_{0}$ be any real number.
\begin{itemize}
\item[a)] The following condition implies condition (\ref{eq:extension_definition_coherence}):
\begin{eqnarray}
\label{eq:extension_condition_one}
\lambda_{0}\leq\inf(-I|B)\Rightarrow\lambda_{0}\leq\inf(-E|B).
\end{eqnarray}
\item[b)] Condition (\ref{eq:extension_definition_coherence}) implies that
\begin{eqnarray}
\label{eq:extension_condition_two}
\lambda_{0}\leq\inf(-I|B)\Rightarrow\lambda_{0}\leq\sup(-E|B)
\end{eqnarray}
\end{itemize}
\end{proposition}
\proof
\begin{itemize}
\item[a)] Let (\ref{eq:extension_condition_one}) hold, and take $\lambda_{0}=\inf(-I|B)$.
Then $\inf(-E|B)\geq\inf(-I|B)=-\sup(I|B)$,
that is $\inf(-E|B)+\sup(I|B)\geq 0$.
We obtain from this $0\leq\sup(\inf(-E|B)+I|B)\leq\sup(I-E|B)$,
which is (\ref{eq:extension_definition_coherence}).
\item[b)] Let (\ref{eq:extension_definition_coherence}) hold.
We obtain, assuming that $\lambda_{0}\leq\inf(-I|B)$ at the third inequality,
$0\leq\sup(I-E|B)\leq\sup(I|B)+\sup(-E|B)=\sup(-E|B)-\inf(-I|B)\leq\sup(-E|B)-\lambda_{0}$.
Hence $\lambda_{0}\leq\sup(-E|B)$, so that (\ref{eq:extension_condition_two}) holds.
\fineproof
\end{itemize}

When $B_{0}=\ldots=B_{n}=\Omega$, i.e. when we consider a bet on unconditional random variables only,
both (\ref{eq:extension_condition_one}) and (\ref{eq:extension_condition_two})
reduce to (\ref{eq:rewrite_alternative_coherence_ter}).
As a by-product,
we reobtain the known result that Definitions \ref{def:coherent_lower_prevision}
and \ref{def:equivalent_definition_coherence} are equivalent.

A comparison of conditions (\ref{eq:rewrite_alternative_coherence_bis}),
(\ref{eq:extension_condition_one}) and (\ref{eq:extension_condition_two}) reveals that
\emph{the expense $E$ is random in a conditional environment}:
it depends on the outcomes of $B_{0},\ldots,B_{n}$ which (apart from those $B_{i}=\Omega$, if any)
are unknown to the bettor at the betting time.
This fact appears to be the real difficulty in trying to extend Definition \ref{def:equivalent_definition_coherence}
to a conditional form:
we actually get two versions, (\ref{eq:extension_condition_one}) and (\ref{eq:extension_condition_two}),
with weaker properties.
Condition (\ref{eq:extension_condition_two}) is potentially useful to disprove W-coherence:
if it does not hold for some bet, the given $\LP(\cdot|\cdot)$ is not W-coherent.
Condition (\ref{eq:extension_condition_one}) is sufficient for W-coherence,
when holding for \emph{any} bet.
A condition slightly simpler than (\ref{eq:extension_condition_one}) may be used
for the same purpose under an additional constraint, as follows
\begin{proposition}
Consider a bet in Definition \ref{def:conditional coherent lower prevision} such that $\wedge_{i=0}^{n}B_{i}\neq\emptyset$.
The following condition implies condition (\ref{eq:extension_definition_coherence}):
\begin{eqnarray}
\label{eq:condition_simpler}
s_{0}\LP(X_{0}|B_{0})-\sum_{i=1}^{n}s_{i}\LP(X_{i}|B_{i})\geq\sup(-I|B).
\end{eqnarray}
\end{proposition}
\proof
We equivalently prove that if (\ref{eq:extension_definition_coherence}) does not hold,
then $p^{*}=s_{0}\LP(X_{0}|B_{0})-\sum_{i=1}^{n}s_{i}\LP(X_{i}|B_{i})<\sup(-I|B)$.
Noting for this that $\wedge_{i=0}^{n}B_{i}\neq\emptyset$ ensures that $p^{*}$ is a possible
value for $-E|B$ and hence $p^{*}\in[\inf(-E|B),\sup(-E|B)]$,
we get $0>\sup(I-E|B)=\sup(-E-(-I)|B)\geq\sup(-E|B)-\sup(-I|B)\geq p^{*}-\sup(-I|B)$,
from which $p^{*}<\sup(-I|B)$ follows.
\fineproof

To ensure W-coherence using (\ref{eq:condition_simpler}), it is necessary that $\wedge_{i=0}^{n}B_{i}\neq\emptyset$,
for any bet.
A relevant special case which obeys this constraint is that of the conditioning events in $\dset$
forming a monotone (or nested) family, i.e. they can be totally ordered by implication
(or inclusion, in the set-theoretic approach).

Summing up, it does not seem possible to generalise Definition \ref{def:equivalent_definition_coherence}
while conditioning.
This should be ascribed to the nature of the term representing the `expense' in the gain decomposition,
which is generally random outside the unconditional framework.

\subsection{Alternative Concepts of Coherence}
\label{sec:alternative_concepts_coherence}
A further issue is that a number of different generalisations of coherence
(Definition \ref{def:coherent_lower_prevision} or equivalent)
to a conditional framework have been proposed in \cite{wa91}:
how do they relate to W-coherence?
We discuss some basic facts about this relationship in this section and the next one.
A further discussion of Walley's criticism on Williams coherence needs some
preliminaries on the concept of avoiding uniform loss, and is therefore presented in Section \ref{sec:aul_condition}.

The coherence concepts defined in Walley's book \cite{wa91} include:
\emph{separate coherence}, which is the first coherence notion in a conditional framework,
presented in Section 6.2.2,
\emph{coherence with unconditional previsions} (Section 6.3.2),
which is generalised to \emph{coherence} in Section 7.1.4 (b),
and \emph{weak coherence}, defined in Section 7.1.4 (a).
Coherence as defined in Section 7.1.4 (b) is the prevailing concept in \cite{wa91},
and will be referred to as \emph{Walley-coherence} here.

None of these concepts \emph{is structure-free}: a common feature is that the conditioning events must belong
to some partition and every (non-impossible) event $B$ in the partition is a conditioning event for some $X|B\in\dset$.
Precisely, just one partition is employed in the case of separate coherence (cf. Definition \ref{def:separate_coherence}),
a \emph{finite} number of partitions are used with Walley-coherence or weak coherence,
two partitions (one of which is the trivial partition $\bset_{0}=\{\Omega\}$,
i.e. it corresponds to unconditional random variables) in the case of coherence with unconditional previsions.
The reason for this kind of constraint lies in Walley's requirement for conglomerability,
a concept discussed in the next section which is itself not structure-free.
There are also other constraints, see e.g. Section 6.3.1,
which are less fundamental,
in the sense that several of them are made to simplify the theory but could be removed;
\cite{mir05} is a paper in this direction.

It ensues that the discussion of, say, Walley-coherence of assignments on relatively simple domains,
like $\dset=\{X_{1}|B_{1}, X_{2}|B_{2}, X_{3}|(B_{1}\wedge(B_{2}\vee B_{3}))\}$, cannot be performed unless these domains
are embedded in larger ones, satisfying the constraints in \cite{wa91}
(this operation could be not simple, it may require some extension theorem which is not always available for Walley-coherence).

Because of these features, Walley's notions of coherence are not always comparable with W-coherence:
there are domains where these notions are not defined,
while W-coherence always is.
When making comparisons,
we must consider W-coherence only on those domains $\dset$ which obey the constraints of the coherence notion
it is compared with.
When this is done, W-coherence is \emph{equivalent} to:
\begin{itemize}
\item[a)] separate coherence (this is proven in Proposition \ref{pro:equivalence_separate_Williams_coherence} below);
\item[b)] Walley-coherence, with the extra assumption that all partitions $\bset_{i}$ of conditioning events
in that definition are finite
(this equivalence is stated without proof in \cite{wa91});
without this assumption, W-coherence is more general than Walley-coherence.
\end{itemize}
As for coherence with unconditional previsions, it is a special case of Walley-coherence.
Concerning weak coherence, it is implied by Walley-coherence
but its importance seems essentially instrumental in the theory in \cite{wa91}.
Useful results for interpreting the conceptual difference between weak coherence and Walley-coherence
were recently given in \cite{mir08}.

Separate coherence has an important role in \cite{wa91},
as it is a prerequisite for the other kinds of coherence.
We are going to prove now its equivalence with W-coherence.
We first state a preliminary result,
which is of some interest in itself, as it simplifies checking W-coherence of $\LP:\dset\rightarrow\reals$
if the conditioning events of all $X|B\in\dset$ have a special separation structure.
\begin{proposition}
\label{pro:partition_coherence}
Given $\LP:\dset\rightarrow\reals$, let $\cset$ be a partition and suppose that, for any $X|B\in\dset$,
$B$ implies some event in $\cset$.
Define $\forall C\in\cset$, $\dset_{C}=\{X|B\in\dset:B\Rightarrow C\}$.
If $\LP$ is W-coherent on each $\dset_{C}$, then it is W-coherent on $\dset$.
\end{proposition}
\proof
The assumptions imply that $\dset=\bigcup_{C\in\cset}\dset_{C}$,
and that a generic gain $G$ in Definition \ref{def:conditional coherent lower prevision}
may be written as follows,
emphasising that distinct random variables may have the same conditioning event:
$G=\sum_{i=1}^{n}\sum_{j=1}^{n_{i}}s_{ij}B_{i}(X_{ij}-\LP(X_{ij}|B_{i}))-s_{0}B_{0}(X_{0}-\LP(X_{0}|B_{0}))$.

Now take, say, $B_{1}$ and suppose $B_{1}\Rightarrow C_{1}\in\cset$.
Then obviously $\sup(G|B)\geq\sup(G|\bigvee_{B_{i}\Rightarrow C_{1}}B_{i})$,
where $\bigvee_{B_{i}\Rightarrow C_{1}}B_{i}$ ($\neq\emptyset$, at least $B_{1}\Rightarrow C_{1}$)
sums those $B_{i}$ among $B_{0},B_{1},\ldots,B_{n}$ that imply $C_{1}$.
But $G|\bigvee_{B_{i}\Rightarrow C_{1}}B_{i}$ is the conditional gain of a bet on (some) elements of $\dset_{C_{1}}$ only,
because those (and only those) $X_{ij}|B_{i}$ (or possibly $X_{0}|B_{0}$) which are not in $\dset_{C_{1}}$ are filtered out,
when conditioning on $\bigvee_{B_{i}\Rightarrow C_{1}}B_{i}$,
by their indicators $B_{i}$ (or $B_{0}$) which all take value zero.
(For instance, if $\bigvee_{B_{i}\Rightarrow C_{1}}B_{i}=B_{1}\vee B_{3}$, $G|B_{1}\vee B_{3}=\sum_{j=1}^{n_{1}}s_{1j}B_{1}(X_{1j}-\LP(X_{1j}))+
\sum_{j=1}^{n_{3}}s_{3j}B_{3}(X_{3j}-\LP(X_{3j}))|B_{1}\vee B_{3})$.
It follows from W-coherence of $\LP$ on $\dset_{C_{1}}$ that $\sup(G|\bigvee_{B_{i}\Rightarrow C_{1}}B_{i})\geq 0$,
hence also $\sup(G|B)\geq 0$.
\fineproof

\textbf{Remark.}
We may replace `W-coherent' with `dF-coherent' in Proposition \ref{pro:partition_coherence}, getting another true proposition.
This is because the preceding proof relies essentially on the structure of $\dset$.
\fineproof

In the sequel we shall apply Proposition \ref{pro:partition_coherence} in the special case that the events $B$ themselves
form partition $\cset$.
Let now $\bset$ be an \emph{arbitrary} (finite or not) partition of non-impossible events.
\begin{definition}
\label{def:separate_coherence}
The conditional lower previsions $\LP_{B}(X|B)$, defined for any $B\in\bset$ and $X\in\hset(B)$,
where $\hset(B)$ is an arbitrary set of random variables containing $B$, are \emph{separately coherent} iff,
for every $B\in\bset$,
\begin{itemize}
\item[i)] $\LP_{B}(B|B)=1$
\item[ii)] $\forall s_{0},\ldots,s_{n}\geq 0$, $\forall X_{0},\ldots,X_{n}\in\hset(B)$,
defining $G=\sum_{i=1}^{n}s_{i}(X_{i}-\LP_{B}(X_{i}|B))-s_{0}(X_{0}-\LP_{B}(X_{0}|B))$,
it holds that $\sup G\geq 0$.\footnote{This is the definition in \cite{wa91},
after replacing integer stakes with real non-negative ones.}
\end{itemize}
\end{definition}
Define now the conditional lower prevision $\LP$ such that $\LP(X|B)=\LP_{B}(X|B)$, $\forall B\in\bset$,
$\forall X\in\hset(B)$ ($\LP$ is the collection of all $\LP_{B}$).
\begin{proposition}
\label{pro:equivalence_separate_Williams_coherence}
The lower previsions $\LP_{B}$ ($B\in\bset$) in Definition \ref{def:separate_coherence} are separately coherent
iff $\LP$ is W-coherent on $\dset=\cup_{B\in\bset}\dset_{B}$, where $\dset_{B}=\{X|B: X\in\hset(B)\}$.
\end{proposition}
\proof
We prove first that W-coherence implies separate coherence.
If $\LP$ is W-coherent, i) trivially holds.
With regard to ii), it follows from

$$\sup G=\max\{\sup_{B}G,\sup_{\setneg{B}}G\}\geq\sup_{B}G=\sup(G|B)=\sup(BG|B)\geq 0,$$
the last equality holding by (\ref{eq:equivalence_conditional_function}), the inequality by W-coherence.

To prove the converse implication, suppose that separate coherence holds.
Betting on $B$, $X_{0},\ldots,X_{n}\in\hset(B)$, it follows then
$\sup(s(B-\LP(B|B))+\sum_{i=1}^{n}s_{i}(X_{i}-\LP(X_{i}|B))-s_{0}(X_{0}-\LP(X_{0}|B)))=
\sup(s(B-1)+G)=
\linebreak
\max(\sup_{B}(s(B-1)+G),\sup_{\setneg{B}}(s(B-1)+G))\geq 0$.

If we choose $s>\max(\sup_{\setneg{B}}G,0)$,
the last inequality implies $\sup_{B}(s(B-1)+G)\geq 0$,
since then $\sup_{\setneg{B}}(s(B-1)+G)=-s+\sup_{\setneg{B}}G<0$.
Using also (\ref{eq:equivalence_conditional_function}),
$\sup_{B}(s(B-1)+G)=\sup(G|B)=\sup(BG|B)=
\sup(\sum_{i=1}^{n}s_{i}B(X_{i}-\LP(X_{i}|B))-s_{0}B(X_{0}-\LP(X_{0}|B))|B)\geq 0$,
which means, given the arbitrariness of $n$, $X_{0},\ldots,X_{n}$ and $s_{0},\ldots,s_{n}\geq 0$,
that $\LP$ is W-coherent on $\dset_{B}$.
Then W-coherence of $\LP$ on each $\dset_{B}$ implies W-coherence of $\LP$
on $\dset$, because of Proposition \ref{pro:partition_coherence}
(where $\cset$, $\dset_{C}$ are now $\bset$, $\dset_{B}$ respectively).
\fineproof

W-coherence and Walley-coherence are equivalent (cf. b) above)
when the partitions $\bset_{i}$ of conditioning events in Walley-coherence are all finite.
In general, properties of W-coherence involving only finitely many distinct conditioning events
hold for Walley-coherence too (a W-coherent assessment or possibly one of its W-coherent extensions,
cf. Proposition \ref{pro:extension_W_coherent},
may be referred in this case to a finite set of finite partitions $\bset_{i}$).
For instance, several product or sign rules are discussed in \cite{pel07} using W-coherence,
but they hold with Walley-coherence too.
One such rule is that, if $\LP$ is W-coherent on $\dset\supset\{AX|B, A|B, X|A\wedge B\}$
and $\LP(X|A\wedge B)>0$,
then
$\LP(AX|B)\geq\LP(A|B)\cdot\LP(X|A\wedge B)$.

In general, W-coherence has the advantage over Walley-coherence that
it verifies properties A), B), D) in Section \ref{sec:historical_note},
while none of them necessarily holds with Walley-coherence.
Property D) allows also a sensitivity analysis interpretation of W-coherence.
W-coherence is not necessarily conglomerative, while Walley-coherence is.
This is a basic difference, and we comment on it in the next Section \ref{sec:non_conglomerability}.

Last but not least, we note that the notion of conditional random variable (and of conditional event)
is often left at an informal level in the literature,
including \cite{wa91,wil07}. A formal approach to these and other descriptive tools of uncertainty,
only sketched in Section \ref{sec:preliminaries},
 is developed in \cite{cri96,cri06}.

Although the way conditional random variables or events are interpreted
is seemingly not particularly relevant in many matters, a greater formalisation
turns out to be useful with other ones.
For an example, consider Lemma 6.2.4 in \cite{wa91}: this lemma states that, if $BX=BY$ and
the separate coherence conditions i), ii) of Definition \ref{def:separate_coherence} hold
for a lower prevision $\LP(\cdot|B)$, then $\LP(X|B)=\LP(Y|B)$.
The result depends on the interpretation of conditional lower previsions in \cite{wa91},
which does not formally define conditional random variables.
But using the approach outlined in Section \ref{sec:preliminaries}
and in particular (\ref{eq:equivalence_conditional_function})
with $n=2$, $X_{1}=B$, $X_{2}=X$, $f(B,X)=BX$,
and since $B|B$ (the indicator of event $B$ given that $B$ is true) takes value $1$,
we get $BX|B=(B|B)\cdot(X|B)=X|B$,
thus condition $BX=BY$ alone implies $X|B=Y|B$.
Consequently we achieve the more general result that
$\mu(X|B)=\mu(Y|B)$ whatever the uncertainty measure $\mu$ is, not because of coherence
($\mu$ could even be incoherent), but merely because we are evaluating the same thing.

\subsection{The Issue of Non-Conglomerability}
\label{sec:non_conglomerability}
Suppose that an uncertainty measure $\mu$ is given on a domain $\dset$
which includes a random variable $X$ and the conditional random variables $X|B$,
for all $B$ in a given partition $\bset$.
Then $\mu$ is \emph{conglomerable} (with respect to $X$ and $\bset$) iff
\begin{eqnarray}
\label{eq:conglomerability_measure}
\inf_{B\in\bset}\mu(X|B)\leq\mu(X)\leq\sup_{B\in\bset}\mu(X|B)
\end{eqnarray}
while $\mu$ is \emph{non-conglomerable} if (\ref{eq:conglomerability_measure}) does not hold.
In words, (\ref{eq:conglomerability_measure}) requires $\mu(X)$ to belong to the smallest interval containing all conditional evaluations
$\mu(X|B)$.

When $X$ is (the indicator of) an event and $\mu$ is a precise probability $P$, conglomerability may seem an obvious
property at first sight,
and in fact it holds trivially if the partition $\bset$ is \emph{finite}.
When $\bset$ is infinite, the matter is however much more complicated \cite{sch84}.

It was de Finetti who discovered in his $1930$ paper \cite{def30} that dF-coherent probabilities may be non-conglomerable,
presenting two nice examples supporting this seemingly counterintuitive fact.
His examples were forerunning the theory, as Definition \ref{def:coherent_lower_prevision} was not known at those times.
We reconsider now one of such examples,
showing that the probability it uses is actually dF-coherent.

\textbf{Example}
A number is chosen at random from the set $\natset^{+}$ of positive integers.
Define $\omega_{n}=\mbox{`}n\mbox{ is chosen}\mbox{'}$,
and term $\bset_{0}$ the partition of all $\omega_{n}$, $n\in\natset^{+}$.

If $A$ is the event that an odd number is chosen,
clearly $P(A)=\frac{1}{2}$.
Defining $B_{n}=\omega_{2n-1}\vee\omega_{4n-2}\vee\omega_{4n}$, $\forall n\in\natset^{+}$,
$\bset_{1,2}=\{B_{1},\ldots,B_{n},\ldots\}$ is a partition coarser than $\bset_{0}$, and
$P(A|B_{n})=\frac{1}{3}$, $\forall n$
(any $B_{n}$ says that either one odd number or two even ones are selected,
so $B_{1}=\mbox{`}1,2\mbox{ or }4\mbox{ is chosen'}$, etc.).
It ensues that (\ref{eq:conglomerability_measure}) does not hold, and $P$ is non-conglomerable.

The example is easily generalised, as noted in \cite{def30}, replacing $\bset_{1,2}$ with the partition $\bset_{h,k}$
such that each of its events $B_{1}^{\prime},\ldots,B_{n}^{\prime},\ldots$ implies that one out of $h+k$ numbers is chosen,
$h$ numbers being odd, $k$ even.
Then $P(A|B_{n}^{\prime})=\frac{h}{h+k}\neq\frac{1}{2}=P(A)$, if $h\neq k$: $P$ is non-conglomerable.

To prove that $P$ is dF-coherent on $\dset=\{A,A|B_{1}^{\prime},\ldots,A|B_{n}^{\prime},\ldots\}$,
note that all possible gains in Definition \ref{def:coherent_prevision} are of two disjoint types,
according to whether they include (a bet on) $A$ or not.
For those who do not $\sup(G|\bigvee_{j=1}^{r}B^{\prime}_{i_{j}})\geq 0$,
applying the remark following Proposition \ref{pro:partition_coherence}
(here $\dset_{C}=\{A|B^{\prime}_{n}\}$, $P$ is dF-coherent on $\dset_{C}$ since $\frac{h}{h+k}\in[0,1])$.
A generic $G$ including $A$ may be written, in a way shorter but equivalent to that of Definition \ref{def:coherent_prevision},
as $G=s(A-\frac{1}{2})+\sum_{j=1}^{r}s_{j}B^{\prime}_{i_{j}}(A-\frac{h}{h+k})$,
where $s,s_{1},\ldots,s_{r}$ may take any real value.
Among those events $\omega_{n}$ of partition $\bset_{0}$ such that $\omega_{n}\wedge(B^{\prime}_{i_{1}}\vee\ldots\vee B^{\prime}_{i_{r}})=\emptyset$,
there are some implying $A$,
while others imply $\setneg{A}$.
If $\omega_{n}\Rightarrow A$, $G(\omega_{n})=\frac{1}{2}s$,
when $\omega_{n}\Rightarrow\setneg{A}$, $G(\omega_{n})=-\frac{1}{2}s$.
In all cases, $\max G\geq 0$.
\fineproof

As this example shows, there may be instances where quite natural uncertainty evaluations are consistent, but non-conglomerable.
We believe that in principle non-conglomerability should not be ruled out a priori.

The issue of non-conglomerability is a root difference between Williams' and Walley's approaches to conditional coherence.
Williams, following de Finetti, does not require conglomerability.
Thus, for instance, the probability $P$ in the example is a special case of W-coherent prevision.

Walley asks for conglomerability in the consistency concepts,
other than separate coherence,
he develops in a conditional framework.
These concepts should comply with a conglomerative principle (\cite{wa91}, Section 6.3.3);
technically, his consistency notions implement this principle by including terms like $G(X|\bset)=\sum_{B\in\bset}B(X-\LP(X|B))$
in the expressions of the gains\footnote{We shall meet one such term in Section \ref{sec:walley_ASL},
equation (\ref{eq:sup_condition_ASL}).
}.
These terms are well-defined also when $\bset$ is infinite, because the factors $B$ are the indicators of events in a partition $\bset$.
Thus only one of them is non-null, whatever happens,
and hence the summation is always made up of a single term.
Conglomerability implies then various conditions, similar to (\ref{eq:conglomerability_measure})
(\cite{wa91}, Section 6.5).
In the case of Walley-coherent precise previsions, it implies axiom (C14) in \cite{wa91}, Section 6.5.7, i.e.
\begin{eqnarray}
\label{eq:conglomerability_c14_walley}
P(X)\geq\inf_{B\in\bset}P(X|B).
\end{eqnarray}
Actually, it is proven in \cite{wa91} that (\ref{eq:conglomerability_c14_walley}) is equivalent to Walley-coherence
for precise previsions,
under certain structure constraints on $\dset$.
These constraints imply in particular that $(X, X|B\in\dset)\Rightarrow(-X,-X|B\in\dset)$,
a condition ensuring alone that (\ref{eq:conglomerability_c14_walley}) is equivalent to (\ref{eq:conglomerability_measure}),
since $P(-X)\geq\inf_{B\in\bset}P(-X|B)$ iff $P(X)\leq\sup_{B\in\bset}(X|B)$.

In particular, it ensues from this argument that the probability $P$ in the example
(technically, any of its dF-coherent extensions on a set $\dset^{\prime}$
meeting the structure requirements of Walley-coherence)
is not Walley-coherent.

More generally, non-conglomerable dF-coherent conditional previsions are not Walley-coherent
(they do not satisfy conglomerative conditions like (\ref{eq:conglomerability_c14_walley})).
Note that the term `\emph{linear prevision}' in \cite{wa91}
identifies dF-coherent previsions in the unconditional environment (the first five chapters),
but corresponds to those dF-coherent conditional previsions which are conglomerable in a conditional setting
(see Section 6.5.7 in \cite{wa91}).

The issue of conglomerability allows a more in-depth explanation of the differences between W-coherence and Walley-coherence.
We pinpoint the following items:
\begin{itemize}
\item[a)]
If we wish that an uncertainty measure $\mu$ is conglomerable,
some constraints must be imposed on its domain $\dset$,
as appears already from (\ref{eq:conglomerability_measure}):
if $X|B\in\dset$, then it must hold that $X|B^{\prime}\in\dset\ \forall B^{\prime}$ in some partition including $B$.
In particular, this or analogous constraints seem unavoidable in Walley-coherence,
which is necessarily not structure-free.
\item[b)]
Walley's approach may be interpreted as a thorough investigation of conglomerable imprecise previsions.
It can be adopted, if one feels that imposing conglomerability does not rule out some significant models in
the specific uncertain situation being investigated.
\item[c)]
Conglomerable imprecise previsions have some additional properties,
ensuing from inequalities like (\ref{eq:conglomerability_measure}), (\ref{eq:conglomerability_c14_walley}),
which are helpful in several derivations and problems.
The disadvantage is that they do not always ensure that the envelope theorem holds,
or that there exists a conglomerable natural extension.
\end{itemize}

\section{Beyond Williams Coherence}
\label{sec:beyond_Williams_coherence}
We explore in this section how W-coherence relates to other consistency concepts,
either stronger (Section \ref{sec:Williams_de_Finetti}) or weaker
(Sections \ref{sec:aul_condition}, \ref{sec:centered_convexity}).

\subsection{Between Williams' and de Finetti's Coherence?}
\label{sec:Williams_de_Finetti}
As well-known, coherence for lower previsions (Definition \ref{def:coherent_lower_prevision})
may be obtained formally from dF-coherence (Definition \ref{def:coherent_prevision})
by restricting the number of bets `against' some $X\in\dset$ (unconstrained with dF-coherence) to $m\leq 1$.
The same constraint distinguishes, in a conditional framework,
W-coherence (Definition \ref{def:conditional coherent lower prevision}) from dF-coherence (Definition \ref{def:conditional_precise_prevision}):
with W-coherence we can bet against (at most) one $X_{0}|B_{0}\in\dset$.

A natural question is then: what if we relax this constraint,
for instance asking - to keep the relaxation at its minimum -
that we can bet
`against' at most two $X|B\in\dset$?
Shall we obtain a significant concept of coherence, intermediate between W-coherence and dF-coherence?
The answer is essentially negative, even in an unconditional environment.
For simplicity, we illustrate this case only.
\begin{definition}
\label{def:2_coherence}
$\LP:\dset\rightarrow\reals$ is a
\emph{bi-coherent lower prevision} on $\dset$
iff, for all $n\in\natset$, $\forall\ X_1,\ldots,X_n,Y_1,Y_2 \in
\dset$, $\forall\ s_1,\ldots,s_n,r_1,r_2$ real and non-negative,
defining $G=\sum_{i=1}^{n}s_{i}(X_i-\LP(X_i))-r_1(Y_1-\LP(Y_1))-r_2(Y_2-\LP(Y_2))$,
$\sup G\geq 0$.
\end{definition}
Clearly, any bi-coherent lower prevision satisfies Definition \ref{def:coherent_lower_prevision} as well and
is therefore coherent.
It also avoids sure loss (\cite{wa91}, Section 2.4.4 (a)), like (as well known) any coherent lower prevision.
Further
\begin{proposition}
\label{pro:properties_2_coherence}
Let $\LP:\dset\rightarrow\rset$ be bi-coherent.
\begin{itemize}
\item[a)] If $X, Y, X+Y\in\dset$, then $\LP(X+Y)=\LP(X)+\LP(Y)$.
\item[b)] If $X, \alpha X\in\dset$ ($\alpha\in\rset$), then $\LP(\alpha X)=\alpha\LP(X)$.
\end{itemize}
\end{proposition}
\proof
To prove a), first observe that the coherence of $\LP$ implies $\LP(X+Y)\geq\LP(X)+\LP(Y)$
(\cite{wa91}, Section 2.6.1 (e)).
For the reverse inequality, put $n=1$, $X_{1}=X+Y$, $Y_{1}=X$, $Y_{2}=Y$, $s_{1}=r_{1}=r_{2}=1$
in Definition \ref{def:2_coherence}.

When $\alpha\geq 0$, b) follows from the coherence of $\LP$ (\cite{wa91}, Section 2.6.1 (f)).
Let us suppose $\alpha<0$.
Putting $n=2$, $X_{1}=X$, $X_{2}=\alpha X$, $s_{1}=-\alpha$, $s_{2}=1$, $r_{1}=r_{2}=0$
in the gain in Definition \ref{def:2_coherence}, we get $\LP(\alpha X)\leq\alpha\LP(X)$.
The opposite inequality follows putting $n=0$, $Y_{1}=X$, $Y_{2}=\alpha X$, $r_{1}=-\alpha$, $r_{2}=1$.
\fineproof

Proposition \ref{pro:properties_2_coherence} emphasises that any bi-coherent lower prevision
is linear and homogenous on a large enough domain, i. e. it behaves essentially like a dF-coherent prevision
(cf. (\ref{eq:linearity})).
Actually, any bi-coherent lower prevision is dF-coherent, when the domain on which it is defined
is sufficiently rich, as the following corollary of Proposition \ref{pro:properties_2_coherence}
points out.
\begin{corollary}
\label{cor:2_coherence_on_spaces}
Let $\LP:\dset\rightarrow\rset$ be bi-coherent.
If either $-X\in\dset\ \forall X\in\dset$ or $X+Y\in\dset\ \forall X,Y\in\dset$, then $\LP$ is dF-coherent.
\end{corollary}
\proof
Let $-X\in\dset\ \forall X\in\dset$.
Since $\LP$ avoids sure loss, and $\LP(X)=-\LP(-X)\linebreak\forall X\in\dset$ by Proposition \ref{pro:properties_2_coherence} b),
dF-coherence of $\LP$ follows at once from Theorem 2.8.2 in \cite{wa91}.
Let now $X+Y\in\dset\ \forall X,Y\in\dset$. Since $\LP$ is coherent, $\LP(X)\geq\LP(Y)+\mu$,
$\forall X,Y\in\dset$ such that $X\geq Y+\mu$ (\cite{wa91}, Section 2.6.1 (d)).
Besides, property a) in Proposition \ref{pro:properties_2_coherence} holds.
This implies dF-coherence of $\LP$ by Theorem 2.8.3 in \cite{wa91}.
\fineproof

Nevertheless, a bi-coherent $\LP$ is not necessarily dF-coherent,
when the domain of $\LP$ does not satisfy the closure properties of Corollary \ref{cor:2_coherence_on_spaces},
as illustrated by the following simple example.

\textbf{Example}
Let $\bset=\{\omega_{1},\omega_{2},\omega_{3}\}$ be a partition
and $\LP$ the vacuous coherent lower prevision on $\bset$: $\LP(\omega_{i})=0$ $(i=1,2,3)$.
Actually, $\LP$ is bi-coherent as well.
To show this, we prove that the supremum of any gain in Definition~\ref{def:2_coherence}
is non-negative. It is sufficient to inspect only the gains of the form
$G_{i}=s_{i}(\omega_{i}-\LP(\omega_{i}))-s_{j}(\omega_{j}-\LP(\omega_{j}))-s_{k}(\omega_{k}-\LP(\omega_{k}))=
s_{i}\omega_{i}-s_{j}\omega_{j}-s_{k}\omega_{k}$
($i\neq j\neq k\neq i, i,j,k\in\{1,2,3\}, s_{i},s_{j},s_{k}\geq 0$),
since the non-negativity of the supremum of any other kind of gain in Definition \ref{def:2_coherence} is implied
by the coherence of $\LP$.
Clearly, $G_{i}(\omega_{i})=s_{i}\geq 0\ (i=1,2,3)$, hence $\LP$ is bi-coherent,
although, patently, $\LP$ is not dF-coherent.
\fineproof

We note incidentally that the vacuous lower prevision is not always bi-cohe\-rent, not even on partitions:
if the partition in the example were $\bset^{\prime}=\{\omega_{1},\omega_{2}\}$,
then $\sup G<0$ in Definition \ref{def:2_coherence}
when $G=-\omega_{1}-\omega_{2}=-1$ (i.e. when $n=0$, $Y_{i}=\omega_{i}$, $r_{i}=1$, $i=1,2$).
This also shows that coherence and bi-coherence are not equivalent,
when bi-coherence may differ from dF-coherence.

Those bi-coherent previsions which are not dF-coherent on $\dset$ do not satisfy property B)
in Section \ref{sec:historical_note},
i.e. they do not ensure bi-coherent extensions on any superset $\dset^{\prime}\supset\dset$.
This is shown by the following corollary.
\begin{corollary}
\label{cor:extension_2_coherent}
Let $\LP:\dset\rightarrow\rset$ be bi-coherent and let $\lset$ be any linear space that contains $\dset$.
Then $\LP$ can be bi-coherently extended on $\lset$ if and only if $\LP$ is dF-coherent on $\dset$.
\end{corollary}
\proof
The `if' part follows from the extension theorem for dF-coherent previsions,
the `only if' part from Corollary \ref{cor:2_coherence_on_spaces}
(implying that \emph{any} bi-coherent extension of $\LP$ on $\lset$ is dF-coherent on $\lset$,
hence also on $\dset\subset\lset$).
\fineproof

Thus, for instance, the lower prevision presented in the previous example cannot
be bi-coherently extended to the set of random variables defined on the partition $\bset$.
Corollary \ref{cor:extension_2_coherent} could be further generalised:
there are instances of bi-coherent, but not dF-coherent, lower previsions that
cannot be bi-coherently extended on supersets which are not even linear spaces.
The important message to convey is anyway already clear: bi-coherence is not particularly significant,
because either it coincides with dF-coherence or, when it can differ from dF-coherence,
property B) of Section \ref{sec:historical_note} may not hold, not even in rather common situations.

\subsection{The Condition of Avoiding Uniform Loss}
\label{sec:aul_condition}
In the unconditional case, the most studied consistency condition weaker than coherence
(Definition \ref{def:coherent_lower_prevision}) is that of avoiding sure loss,
obtained formally from Definition \ref{def:coherent_lower_prevision} putting $s_{0}=0$.
With W-coherence, the corresponding weaker notion is the following
\begin{definition}
\label{def:aul_prevision}
$\LP:\dset\rightarrow\reals$
\emph{avoids uniform loss} (AUL)
iff, for all $n\in\natset^{+}$, $\forall\ X_1|B_{1},\ldots,X_n|B_{n} \in
\dset$, $\forall\ s_1,\ldots,s_n$ real and non-negative,
defining $B=\vee_{i=1}^{n}B_{i}$ and $G=\sum_{i=1}^{n}s_{i}B_{i}(X_i-\LP(X_i))$,
it holds that $\sup(G|B)\geq 0$.
\end{definition}
The notion of avoiding uniform loss was used in \cite{wa04}, where other equivalent characterisations are supplied.
When $\LP=\UP=P$, $P$ avoids uniform loss if and only if $P$ is dF-coherent.
Clearly, W-coherence of $\LP$ implies that $\LP$ avoids uniform loss.
The AUL condition is generally too weak, as appears already at the unconditional level (cf. \cite{wa91}, Section 2.5).
A more satisfactory notion is that of centered convexity (cf. Section \ref{sec:centered_convexity}).

In this section we explore the relationship between the AUL condition and a similar notion introduced in \cite{wa91},
and reconsider an example on W-coherence discussed in \cite{wa91} in the light of this.

\subsubsection {Walley's Condition of Avoiding Sure Loss}
\label{sec:walley_ASL}
Given a partition $\bset$ and two arbitrary sets $\hset$, $\kset$ of unconditional random variables,
such that $0\in\hset$, $B\in\hset\ \forall B\in\bset$,
suppose throughout this section that $\dset$ has the following special structure:
$\dset=\kset\cup\bigcup_{B\in\bset}\dset_{B}$, where $\dset_{B}=\{Y|B:Y\in\hset\}$.
\begin{definition}
\label{def:asl_prevision}
Let $\LP:\dset\rightarrow\reals$ be such that
\begin{itemize}
\item[a)] the restriction of $\LP$ on $\kset$ is a(n unconditional) coherent lower prevision;
\item[b)] the restrictions of $\LP$ on each $\dset_{B}$, $B\in\bset$, are separately coherent.
\end{itemize}
Then $\LP$ \emph{avoids sure loss} on $\dset$
iff, for all $m, n\in\natset$, $\forall\ X_1,\ldots,X_m \in
\kset$, $\forall\ Y_1,\ldots,Y_n \in \hset$,$\forall\ s_j\geq 0\ (j=1,\ldots,m)$,  $\forall\ t_i\geq 0\ (i=1,\ldots,n)$,
\begin{eqnarray}
\label{eq:sup_condition_ASL}
\sup(\sum_{j=1}^{m}s_{j}(X_j-\LP(X_j))+\sum_{i=1}^{n}t_{i}\sum_{B\in\ \bset}B(Y_i-\LP(Y_i|B)))\geq 0.
\end{eqnarray}
\end{definition}
\textbf{Discussion.}
Definition \ref{def:asl_prevision} is Walley's Definition 6.3.2 of \emph{avoiding sure loss} in \cite{wa91};
here the following assumptions\footnote{Condition (b) of Section 6.3.1 in \cite{wa91},
i.e. $Y\in\hset\Rightarrow BY\in\hset,\ \forall B\in\bset$, may be replaced by our assumptions on $\dset$,
in particular by $0\in\hset$.
In fact, given any $B,B^{*}\in\bset$, we have that $B^{*}Y|B$ is equal to $Y|B$, by (\ref{eq:equivalence_conditional_function}),
when $B^{*}=B$, while, when $B^{*}\neq B$, $B^{*}Y|B=0|B\ (\in\dset)$.
Therefore, $\LP(B^{*}Y|B)$ is defined $\forall B,B^{*}\in\bset$, which is what ensures condition (b) of Section 6.3.1 in \cite{wa91}.
We did not mention condition (c) of Section 6.3.1 because it is unnecessary in the following derivations.
}
are introduced, without altering Definition 6.3.2:
\begin{itemize}
\item[i)] the non-negative coefficients $s_{j}$, $t_{i}$ are real, but not necessarily integers;
\item[ii)] we do not require $\hset$, $\kset$ to be linear spaces,
unlike condition (a) in \cite{wa91}, Section 6.3.1,
and modified correspondingly Definition 6.3.2, as indicated at the end of Section 6.3.1.
\end{itemize}

An interesting remark is that \emph{Definition \ref{def:asl_prevision} is formally no extension of the condition
of avoiding sure loss for unconditional previsions} (Definition 2.4.4 (a) in \cite{wa91}):
if $\kset=\varnothing$ and $\bset=\{\Omega\}$,
it reduces to the notion of coherence (Definition \ref{def:coherent_lower_prevision}).
This depends on assuming b) in Definition \ref{def:asl_prevision}.
The same remark applies to the concepts of avoiding sure, partial or uniform loss\footnote{Note that
the meaning of the term avoiding uniform loss in \cite{wa91} is different from that used in this paper,
following Definition~\ref{def:aul_prevision}.
}
defined in \cite{wa91},
chapter 7, since separate coherence is a prerequisite for them too.
\fineproof

\begin{proposition}
\label{pro:ASL_implies_AUL}
If $\LP$ avoids sure loss on $\dset$, it avoids uniform loss on $\dset$.
\end{proposition}
\proof
Given the special structure of $\dset$, all gains $G$ in Definition \ref{def:aul_prevision} may be written as follows,
\begin{eqnarray}
\label{eq:gain_AUL}
G=\sum_{j=1}^{m}s_{j}(X_j-\LP(X_j))+\sum_{i=1}^{q}\sum_{r=1}^{k_{i}}t_{ir}B_{i}(Y_{ir}-\LP(Y_{ir}|B_i)),
\end{eqnarray}
with $m, q\geq 0$.
Suppose $\LP$ avoids sure loss, and consider the following (exhaustive) cases.
\begin{itemize}
\item[i)] The second summation in (\ref{eq:gain_AUL}) is zero ($q=0$).
Then $\sup G\geq 0$ follows from Definition \ref{def:asl_prevision}, a).
\item[ii)] The first summation in (\ref{eq:gain_AUL}) is zero ($m=0$).
Separate coherence of $\LP$ on all $\dset_{B}$ (Definition \ref{def:asl_prevision}, b)) implies W-coherence
of $\LP$ on $\bigcup_{B\in\bset}\dset_{B}$ (Proposition \ref{pro:equivalence_separate_Williams_coherence}),
which implies that $\LP$ avoids uniform loss on $\bigcup_{B\in\bset}\dset_{B}$, hence $\sup G|\vee_{i=1}^{q}B_{i}\geq 0$.
\item[iii)] $m\cdot q>0$.
This implies $\sup(G|B)=\sup(G|\Omega)=\sup G$ in Definition \ref{def:aul_prevision}.
We can write $G$ as a gain of the kind (\ref{eq:sup_condition_ASL}),
since $B_{i}(Y_{ir}-\LP(Y_{ir}|B_i))=\sum_{B\in\bset}B(B_{i}Y_{ir}-\LP(B_{i}Y_{ir}|B))$
(we used the fact that $BB_{i}=0$ if $B\neq B_{i}$, and that $B_{i}Y_{ir}|B=B_{i}|B\cdot Y_{ir}|B$;
consequently if $B\neq B_{i}$, $B_{i}Y_{ir}|B=0|B$, and $\LP(B_{i}Y_{ir}|B)=0)$.
Then $G$ in (\ref{eq:gain_AUL}) is a gain of type (\ref{eq:sup_condition_ASL})
from a bet on $X_{1},\ldots,X_{m}$, and on the conditional random variables
$B_{1}Y_{11}|B,\ldots,\linebreak B_{q}Y_{qk_{q}}|B$, $\forall\ B\in\bset$.
Hence $\sup G\geq 0$.
\end{itemize}
In all cases, $G$ satisfies the conditions in Definition \ref{def:aul_prevision}.
\fineproof

Hence, Definition \ref{def:asl_prevision} is stronger than Definition \ref{def:aul_prevision},
when they are comparable.
The key difference is that Definition \ref{def:asl_prevision} can be justified following a conglomerative principle
(cf. \cite{wa91}, Section 6.3.3) while Definition \ref{def:aul_prevision} does not rely on it.
This fact is relevant in explaining some of Walley's remarks on Williams coherence, as we shall now see.

\subsubsection{On the Consistency of Williams Coherence}
\label{sec:consistency_Williams}
A critical remark in \cite{wa91} about Williams coherence is that it does not always satisfy Walley's condition
of avoiding sure loss.

The important fact here is that if an agent adopts W-coherence, her/his reference minimal consistency concept
should be Definition \ref{def:aul_prevision} of avoiding uniform loss, or equivalent.
Referring to Definition \ref{def:asl_prevision} of avoiding sure loss would determine a kind of inconsistency:
the agent requires (with the condition of avoiding sure loss)
and does not require (with W-coherence) conglomerability at the same time.

Keeping the concept of avoiding uniform loss as a reference,
the criticism to W-coherence outlined in some examples in \cite{wa91} does not apply.
We discuss here one such example (\cite{wa91}, Section 6.6.6).

Let $\bset$ be a denumerable partition whose elements are indexed in the set $\intset-\{0\}$ of non-zero integers and
call $\omega_{z}$ the generic element in $\bset$.
Define two dF-coherent precise probabilities $P^{+}$ and $P^{-}$ on $\bset$, as follows.
$P^{+}(\omega_{z})=2^{-z}$ if $z>0$, $P^{+}(\omega_{z})=0$, $\forall\ z<0$,
while $P^{-}(\omega_{z})=0$, $\forall z$.
Extend $P^{+}$, $P^{-}$ on $B=\bigvee_{\{z<0\}}\omega_{z}$:
clearly $P^{+}(B)=0$ ($P^{+}$ is $\sigma$-additive),
while the extension of $P^{-}$ is not unique, and we may dF-coherently choose $P^{-}(B)=1$.
The extensions on $A=\setneg{B}=\bigvee_{\{z>0\}}\omega_{z}$ are then $P^{-}(A)=0$, $P^{+}(A)=1$.

Define now $P=\frac{P^{+}+P^{-}}{2}$. Since mixtures of dF-coherent probabilities are dF-coherent,
$P$ is dF-coherent.
Let $n\in\nset^{+}$, and define $B_{n}=\omega_{-n}\vee\omega_{n}$.
Because $P(B_{n})=P(\omega_{n})>0$, the extension of $P$ on $\omega_n|B_{n}$ is uniquely determined by Bayes' rule,
and $P(\omega_n|B_{n})=1$.
Similarly, $P(A|B_{n})=1$, while $P(A)=\frac{1}{2}$.
Then $P$ \emph{is} a dF-coherent conditional probability on $\dset=\bset\cup\{B, A, B_{n}, \omega_n|B_{n}, A|B_{n}\}$:
this follows from the fact that coherent (conditional or not) probabilities can be
dF-coherently extended on any event \cite{cri06,def74,hol85},
and that the extension of $P$ on $\omega_n|B_{n}$ and $A|B_{n}$ is unique.
DF-coherence of $P$ on $\dset$ is equivalent to its avoiding uniform loss on $\dset$,
when viewing $P$ as a special imprecise prevision \cite{wa04}.
Thus $P$ does not incur uniform loss, but it is shown in \cite{wa91} that it incurs sure loss
(in the sense of Definition \ref{def:asl_prevision}).
This is because $P$ is non-conglomerable, and in fact it does not obey the conglomerability axiom
(\ref{eq:conglomerability_c14_walley}).

Similar conclusions hold for other examples in \cite{wa91}:
inconsistencies arise \emph{only} when conglomerability axioms are used in a hybrid way.
Thus the very question in choosing between W-coherence and Walley-coherence (when they do not coincide)
seems to be a problem of imposing or not conglomerability.

\subsection{Centered Convexity}
\label{sec:centered_convexity}
While modifications of the definition of W-coherence towards some notions intermediate between it and dF-coherence
seem to yield no really significant results,
the notion of centered convexity is intermediate between that of avoiding uniform loss and W-coherence and has
interesting properties.
\begin{definition}
\label{def:conditional convex lower prevision}
$\LP:\dset\rightarrow\reals$ is a \emph{convex conditional lower prevision} on $\dset$ iff
$\forall n\in\natset^{+}$,
$\forall X_0|B_0,\ldots,X_n|B_n\in\dset$, $\forall s_1,\ldots,s_n\geq 0:
\sum_{i=1}^{n}s_i=1$, defining $G=\sum_{i=1}^{n}s_{i}B_{i}(X_{i}-\LP(X_{i}|B_{i}))-B_{0}(X_{0}-\LP(X_{0}|B_{0}))$,
$\sup\{G|\vee_{i=0}^{n}B_{i}\}\geq 0$.
Further, $\LP$ is \emph{centered} if besides $0|B\in\dset$ and $\LP(0|B)=0,\forall X|B\in\dset$.
\end{definition}
The theory of centered convex previsions was developed in \cite{pel03bis,pel05,pel05bis}, generalising under many respects
the theory of W-coherence.
These previsions satisfy the properties A), B) and C) from Section 2.1,
and operationally correspond to the important notion of convex risk measure.

Property D) in Section \ref{sec:historical_note} is the only one, among those stressed in this paper,
where W-coherence still has a definite advantage over centered convexity, at the current state of art.
In the rest of this section, we give some explanation of this fact.
The material is derived from \cite{pel05bis},
where the interested reader may find more details.
We present here the simplest envelope theorem,
whose proof requires preliminarily the following
\begin{proposition}
\label{pro:envelope_convexity}
Let $\pset$ be a set of convex conditional lower previsions defined on $\dset$.
If $\LP(X|B)=\inf_{\LQ\in\pset}\left\{\LQ(X|B)\right\}$ is finite $\forall X|B\in\dset$,
$\LP$ is convex on $\dset$.
\end{proposition}
Proposition \ref{pro:envelope_convexity} generalises to convex conditional lower previsions a statement already
established for coherent \cite{wa91} or convex unconditional \cite{pel03bis} lower previsions.\footnote{\label{foo:conceptual_difference}There
is a conceptual difference with coherence: since convexity does not imply $\LQ(X|B)\geq\inf(X|B)$ (internality),
the finiteness condition of the infimum must be required in Proposition \ref{pro:envelope_convexity}.
Internality holds when the convex previsions are centered.
This fact exemplifies that convexity without centering may be a rather weak consistency requirement.
}
The proof is similar to those in \cite{pel03bis, wa91} and is omitted.

\textbf{Notation}
Given $\dset$, let $\eset=\{B:\exists X|B\in\dset\}$.
\fineproof

\begin{theorem}
\label{thm:envelope_theorem_without_structure}
\emph{(Envelope Theorem)} Let $\pset$ be a set of dF-coherent precise previsions on $\dset\cup\eset$
such that $\forall P\in\pset$, $P(B)>0\ \forall B\in\eset$,
and let $\alpha:\pset\rightarrow\rset$ be a real function.
Then
\begin{eqnarray}
\label{eqn:envelope_without_structure}
\LP(X|B)=\inf_{P\in\pset}\{P(X|B)+\frac{\alpha(P)}{P(B)}\}\ \forall X|B\in\dset
\end{eqnarray}
is a convex conditional lower prevision on $\dset$, whenever the infimum in (\ref{eqn:envelope_without_structure})
is finite.
Further, $\LP$ is centered iff $\inf_{P\in\pset}\{\frac{\alpha(P)}{P(B)}\}=0,\ \forall B\in\eset$.
\end{theorem}
\proof
We prove that $\forall P\in\pset$, $\forall \alpha\in\rset$, $\LP_{\alpha}=P(X|B)+\frac{\alpha}{P(B)}$ is convex.
The main thesis of the theorem then follows from Proposition \ref{pro:envelope_convexity}.

To prove that $\LP_{\alpha}$ is a convex conditional lower prevision,
we show that a generic $G$ in Definition \ref{def:conditional convex lower prevision}
may be referred to $P$,
after substituting $\LP_{\alpha}(X|B)$ with $P(X|B)+\frac{\alpha}{P(B)}$,
and hence its supremum is non-negative because $P$ is dF-coherent.
In fact, let $X_0|B_0,\ldots,X_n|B_n\in\dset$, $s_1,\ldots,s_n\geq 0$
such that $\sum_{i=1}^{n}s_i=1$.
Then $G$ can be written as

$G=\sum_{i=1}^{n}s_{i}B_{i}(X_{i}-P(X_{i}|B_{i})-\alpha/P(B_{i}))-B_{0}(X_{0}-P(X_{0}|B_{0})-\alpha/P(B_{0}))=$

$\sum_{i=1}^{n}s_{i}B_{i}(X_{i}-P(X_{i}|B_{i}))+\sum_{i=1}^{n}s_{i}(B_{i}\vee B_{0})(Z_{i}-P(Z_{i}|B_{i}\vee B_{0}))
-B_{0}(X_{0}-P(X_{0}|B_{0}))$,
where $Z_{i}=\alpha(B_{0}/P(B_{0})-B_{i}/P(B_{i}))$ and
$P(Z_{i}|B_{i}\vee B_{0}))=
\alpha(P(B_{0}|B_{i}\vee B_{0})/P(B_{0})-P(B_{i}|B_{i}\vee B_{0})/P(B_{i}))=
\alpha(1/P(B_{i}\vee B_{0})-1/P(B_{i}\vee B_{0}))=0$
is, by (\ref{eqn:linearity_coherence}),
the only coherent extension of $P$ on $Z_{i}|B_{i}\vee B_{0}$, $i=1,\ldots,n$.
In terms of $P$, the gain $G$ is still conditioned on $B$,
because $B$ is also the logical sum of the new conditioning events:
$B=\bigvee_{i=1}^{n}B_{i}\vee\bigvee_{i=1}^{n}(B_{i}\vee B_{0})\vee B_{0}$.
It follows $\sup G|B\geq 0$ by dF-coherence of $P$.

The proof of the second part of the proposition follows at once from noting that
when $X|B=0|B$ (\ref{eqn:envelope_without_structure}) reduces to
$\LP(0|B)=\inf_{P\in\pset}\{\frac{\alpha(P)}{P(B)}\}$.
\fineproof

Theorem \ref{thm:envelope_theorem_without_structure} is not a characterisation theorem,
and cannot obviously be applied to arbitrary $\dset$ and $\eset$.
One reason for presenting it is that it supplies us with a way of assessing centered convex previsions
in the particular, but important case that $P(B)>0$, $\forall\ B\in\bset$, $\forall\ P\in\pset$.

Another motivation is that it informs us, through (\ref{eqn:envelope_without_structure}),
about the type of functions upon which the infimum is performed.
Convexity requires adding a term $\phi_{P}(B)$ to any dF-coherent prevision $P(X|B)$.
This term is equal to $\frac{\alpha(P)}{P(B)}$ in Theorem \ref{thm:envelope_theorem_without_structure}.
If $\LP$ is unconditional, it reduces to $\alpha(P)$,
if it is W-coherent, $\phi_{P}\equiv 0$,
and we come to the familiar envelope theorems in \cite{wa91, wil07}.

An envelope theorem which characterises centered convexity is given in \cite{pel05bis}, Theorem 8.
We do not report it here,
but stress the fact that its practical use is considerably less immediate than the envelope theorem for W-coherence.
In fact, the set on which the infimum is performed \emph{depends on} $X|B$ in this theorem.
Also the function $\phi_{P}(B)$ has a more complex structure,
which is influenced by the ordering of zero probabilities, for each $P\in\pset$,
among the possible conditioning events.
Seemingly, it is technically possible to circumvent such difficulties with W-coherence
because the function $\phi_{P}(B)$ may be set identically equal to zero there.

Thus W-coherence remains so far the most general concept for which D) in Section \ref{sec:historical_note}
has a general practical as well as theoretical significance among those discussed in this paper.

\section{Conclusions}
\label{sec:conclusions}
We summarise our conclusions about the role of Williams coherence with the help of Table \ref{ta:puncond},
where consistency concepts for precise (first) and imprecise previsions (then) are listed in order of increasing generality.
\renewcommand{\arraystretch}{1.25}
\begin{table}[ht]
\caption{Some consistency concepts for precise and   imprecise previsions}
\label{ta:puncond}
\renewcommand{\arraystretch}{2.50}
\small
\begin{tabular}{|c|c|c|c|c|}
  \hline
  & {\begin{minipage} {2.5cm} \center{\textbf{Type of Prevision}} \end{minipage}} & \begin{minipage} {1.9cm} \center{\textbf{A) Structure Free}} \end{minipage}&
    \begin{minipage} {2.0cm} \center{\textbf{B) Extension Theorem}} \end{minipage} & \begin{minipage}{3.3cm} \center{\textbf{D) Envelope Theorem Characterisation}} \end{minipage}  \\
  \hline
  \begin{minipage} {2.0cm} \center{de Finetti - coherence} \end{minipage} & \begin{minipage} {2.5cm} \center{Precise, unconditional} \end{minipage} & Yes & Yes & Does not apply \\
  \hline
  \begin{minipage} {2.0cm} \center{de Finetti - coherence} \end{minipage} & \begin{minipage} {2.5cm} \center{Precise, conditional} \end{minipage} & \multicolumn{2}{c|}{Yes, in later studies} & Does not apply \\
  \hline
  \begin{minipage} {2.0cm} \center{Coherence} \end{minipage} & \begin{minipage} {2.5cm} \center{Lower, unconditional} \end{minipage} & Yes & Yes & Yes \\
  \hline
  \begin{minipage} {2.0cm} \center{Walley-coherence} \end{minipage} & \begin{minipage} {2.5cm} \center{Lower, conditional} \end{minipage} & No & Not always & Not always \\
  \hline
  \begin{minipage} {2.0cm} \center{W-coherence} \end{minipage} & \begin{minipage} {2.5cm} \center{Lower, conditional} \end{minipage} & Yes & Yes & Yes \\
  \hline
  \begin{minipage} {2.0cm} \center{Centered convexity} \end{minipage} & \begin{minipage} {2.5cm} \center{Lower, conditional} \end{minipage} & Yes & Yes & \begin{minipage} {2.5cm} \center{Yes (with operational constraints)} \end{minipage}\\
  \hline
 \end{tabular}
\normalsize
\end{table}
%
Undoubtedly, a strong motivation for adopting the variant of Williams coherence called W-coherence in this paper is its generality:
it meets all the properties we listed in Section \ref{sec:historical_note},
a feature shared by coherence for unconditional lower previsions and the root concept of dF-coherence.
Even the notion beyond W-coherence, i.e. centered convexity, while being more general (but weaker)
under many respects, fails to ensure a general envelope theorem of comparable ease of use.
If we restrict our attention to W-coherence versus Walley-coherence,
we may conclude that whenever they are not equivalent (if they are we may adopt either one)
the choice depends essentially on our willingness to accept some conglomerative axiom,
and some at a large extent consequent domain constraints (acceptance of both items results in preferring Walley-coherence).
Given that W-coherence is more general than Walley-coherence,
we may even use W-coherence in principle, and Walley-coherence under specific circumstances,
for instance when studying stochastic processes.
This case copes well with the domain constraints of Walley-coherence,
when we are interested in lower previsions like
$\LP(X_{n}|\wedge_{i=1}^{n-1}(X_{i}=x_{i}))$,
where $x_{i}$ is a generic value for the random variable $X_{i}$.
In fact, the events $\wedge_{i=1}^{n-1}(X_{i}=x_{i})$
form a partition $\bset_{n-1}$, for a given $n$ and by varying $x_{1},\ldots,x_{n-1}$ in all (jointly) possible ways.

Similarly, new information in statistical inference may commonly arise from a partition of possible hypotheses.
Again, this is a favourable situation to apply Walley-coherence, as for its domain constraints,
and is in fact largely discussed in \cite{wa91}.
It has also to be noted that cases where Walley-coherence ensures the existence of a (conglomerable) extension
are pointed out in \cite{wa91},
and that they are of a certain generality.
In other words, the `not always' at the crossing of Walley-coherence and property B) in Table \ref{ta:puncond} should be graded.

More generally, the theory of imprecise probabilities shows that there are often many alternatives for generalising
familiar concepts (for instance, independence) from theories of precise probabilities or previsions,
and that frequently there is no way to keep all the properties of the special precise probability case.
Under these circumstances, we might want to employ different concepts of conditional consistency,
to preserve obtaining certain aims.
A presentation of these conflicting instances is given in \cite{wei03},
where some alternative notions of imprecise conditional probability are presented.
A further investigation of the consistency concepts in the conditional environment should include
also these aspects, as well as other ideas developed in the literature.
In particular, the game-theoretic approach in \cite{sha01,sha03} was recently related to Walley's \cite{dec07},
and this could simplify the potential future work of relating it with Williams' approach too.

\section*{Appendix. Proof of Proposition \ref{pro:extension_W_coherent}.}
We preliminarily recall a characterisation theorem,
holding for W-coherent conditional lower previsions
defined on a structured domain $\dset^{*}$ \cite{wil07}.
\begin{theorem}
\label{thm:conditional coherent prevision on linear space}
Let $\xset$ be a linear space of bounded random variables, $\eset\subset\xset$ the set of all indicator functions
of events in $\xset$.
Let also $1\in\eset$ and $BX\in\xset$, $\forall B\in\eset$, $\forall X\in\xset$.
Define $\eset^{\emptyset}=\eset-\{\emptyset\}$, $\dset^{*}=\{X|B:X\in\xset, B\in\eset^{\emptyset}\}$.
$\LP:\dset^{*}\rightarrow\reals$ is a W-coherent conditional lower prevision if and only if:
\begin{itemize}
\item[A1)] $\LP(X|B)\geq\inf\{X|B\}, \forall X|B\in\dset^{*}$
\item[A2)] $\LP(kX|B)=k\LP(X|B), \forall X|B\in\dset^{*}$, $\forall k\geq 0$
\item[A3)] $\LP(X+Y|B)\geq\LP(X|B)+\LP(Y|B),
\forall X|B,Y|B\in\dset^{*}$
\item[A4)] $\LP(A(X-\LP(X|A\wedge B))|B)=0, \forall X\in\xset$$, \forall A, B\in\eset^{\emptyset}:
A\wedge B\neq\emptyset$.
\end{itemize}
\end{theorem}
As in the unconditional case \cite{wa91},
the concept of \emph{natural extension} plays a fundamental role in extending $\LP$.
\begin{definition}
\label{def:conditional natural extension}
Let $\LP:\dset\rightarrow\reals$ be a conditional lower prevision, $X|B$
an arbitrary bounded conditional random variable. Define
$g_i=s_{i}B_{i}(X_{i}-\LP(X_{i}|B_{i}))$,
$L(X|B)=\{\alpha:\sup\{\sum_{i=1}^{n}g_{i}-B(X-\alpha)|\bigvee_{i=1}^{n}B_{i}\vee
B\}<0, \mbox{ for some } n\geq 0, X_{i}|B_{i}\in\dset, s_{i}\geq
0\}$. The \emph{natural extension} of $\LP$ to $X|B$ is
$\LE(X|B)=\sup L(X|B)$.
\end{definition}
It is easily seen that $L(X|B)=]-\infty,\LE(X|B)[$, a fact which will be used later.
Moreover, the natural extension proves to be bounded from above,
when $\LP$ is W-coherent.
\begin{proposition}
Let $\LP:\dset\rightarrow\reals$ be a W-coherent conditional lower prevision.
Then $\LE(X|B)\leq\sup\{X|B\}\ \forall X|B$.
\end{proposition}
\proof
Let $c=\sup\{X|B\}$, $n\geq 0$, $X_{i}|B_{i}\in\dset$, $s_{i}\geq 0\ (i=1,\ldots,n)$.
Since $B(X-c)\leq 0$, using also W-coherence of $\LP$ in the last inequality,
$\sup\{\sum_{i=1}^{n}g_{i}-B(X-c)|\bigvee_{i=1}^{n}B_{i}\vee B\}\geq
\sup\{\sum_{i=1}^{n}g_{i}|\bigvee_{i=1}^{n}B_{i}\}\geq 0$.
This implies $c\notin L(X|B)=]-\infty,\LE(X|B)[$.
\fineproof

\begin{theorem}
\label{thm:extension_conditional_lower_prevision}
Let $\dset^{*}$ be defined as in Theorem \ref{thm:conditional coherent prevision on linear space},
$\dset\subset\dset^{*}$ and $\LP:\dset\rightarrow\reals$ W-coherent.
Then $\LE$ is a W-coherent conditional lower prevision on $\dset^{*}$
and $\LE(X|B)=\LP(X|B)\ \forall X|B\in\dset$.
\end{theorem}
\proof
To prove W-coherence of $\LE$, we show that it satisfies properties A1), A2), A3), A4) in
Theorem \ref{thm:conditional coherent prevision on linear space}.

As for A1), note that
$\sup\{-B(X-\alpha)|B\}<\sup\{-B(X-\inf\{X|B\})|B\}\leq 0$,
$\forall\ X|B\in\dset^{*}$, $\forall\ \alpha<\inf\{X|B\}$.
This implies $\LE(X|B)\geq \inf\{X|B\}$.

As for A2), let $k>0$ (the case $k=0$ is trivial), $\alpha\in L(X|B)$, $n\geq 0$, $X_{i}|B_{i}\in\dset$, $s_{i}\geq 0\ (i=1,\ldots,n)$,
$W_{1}=\sum_{i=1}^{n}g_{i}-B(X-\alpha)$ as in Definition \ref{def:conditional natural extension}.
Then, $\sup\{\sum_{i=1}^{n}k s_{i}(X_{i}-\LP(X_{i}|B_{i}))-B(kX-k\alpha)|\bigvee_{i=1}^{n}B_{i}\vee B\}=
k\sup\{W_{1}|\bigvee_{i=1}^{n}B_{i}\vee B\}<0$.
This implies $k\alpha\in L(kX|B)\ \forall k>0$, $\forall\alpha\in L(X|B)$.
Hence $\LE(kX|B)\geq k\LE(X|B)$.
The proof of the reverse inequality is similar.

To prove A3), let $Y|B\in\dset^{*}$, $\beta\in L(Y|B)$, $m\geq 0$, $Y_{j}|C_{j}\in\dset$, $t_{j}\geq 0\ (j=1,\ldots,m)$,
$h_{j}=t_{j}C_{j}(Y_{j}-\LP(Y_{j}|C_{J}))$ such that,
defining $W_{2}=\sum_{j=1}^{m}h_{j}-B(Y-\beta)$, $\sup \{W_{2}|\bigvee_{j=1}^{m}C_{j}\vee B\}<0$.
Preliminarily, write $H=\bigvee_{i=1}^{n}B_{i}\vee\bigvee_{j=1}^{m}C_{j}\vee B$ as the sum of four disjoint events as follows:
$H=
B\vee
[\bigvee_{i=1}^{n}B_{i}\wedge\setneg{(\bigvee_{j=1}^{m}C_{j})}\wedge \setneg{B}]\vee
[\setneg{(\bigvee_{i=1}^{n}B_{i})}\wedge\bigvee_{j=1}^{m}C_{j}\wedge \setneg{B}]\vee
[\bigvee_{i=1}^{n}B_{i}\wedge\bigvee_{j=1}^{m}C_{j}\wedge \setneg{B}]$.
Observe also that $\sup W_{1}$, $\sup W_{2}$ are both non-positive, but never simultaneously null,
conditional on each of the four events.
This implies $\sup\{W_{1}+W_{2}|H\}=
\sup\{\sum_{i=1}^{n}g_{i}+\sum_{j=1}^{m}h_{j}-B(X+Y-(\alpha+\beta))|H\}<0$.
Hence $\alpha+\beta\in L(X+Y|B)\ \forall \alpha\in L(X|B), \forall\beta\in L(Y|B)$ and $\LE(X+Y|B)\geq\LE(X|B)+\LE(Y|B)$ follows.

As for A4), let $X|A\wedge B\in\dset^{*}$, $W=A(X-\LE(X|A\wedge B))$.
To prove that $\LE(W|B)=\sup L(W|B)=0$, we show that $L(W|B)=]-\infty,0[$.
Given $\delta>0$, it ensues from the definition of $\LE(X|A\wedge B)$  that $\exists\ n\geq 0$, $X_{i}|B_{i}\in\dset$,
$s_{i}\geq 0\ (i=1,\ldots,n)$ such that, defining
$G=\sum_{i=1}^{n}s_{i}B_{i}(X_{i}-\LP(X_{i}|B_{i}))$ and
$Z_{1}=G-AB(X-\LE(X|A\wedge B)+\delta)$,
$\sup \{Z_{1}|\bigvee_{i=1}^{n}B_{i}\vee(A\wedge B)\}<0$.
Hence $Z_{2}=G-B(W+\delta)=Z_{1}-B\setneg{A}\delta$ $(\leq Z_{1})$ is such that
$\sup \{Z_{2}|\bigvee_{i=1}^{n}B_{i}\vee B\}=\max\{\sup \{Z_{2}|\bigvee_{i=1}^{n}B_{i}\vee(A\wedge B)\},
\sup \{Z_{2}|\setneg{(\bigvee_{i=1}^{n}B_{i})}\wedge\setneg{A}\wedge B\}\}\leq
\max\{\sup \{Z_{1}|\linebreak\bigvee_{i=1}^{n}B_{i}\vee(A\wedge B)\},-\delta\}<0$
(omit the second argument in the maxima if $\setneg{(\bigvee_{i=1}^{n}B_{i})}\wedge\setneg{A}\wedge B=\emptyset$).
This implies $-\delta\in L(W|B),\ \forall \delta>0$, hence $\sup L(W|B)\geq 0$.
But $\sup L(W|B)=0$, because $0\notin L(W|B)$: by contradiction, assuming $0\in L(W|B)$ would imply, as can be easily seen,
$\LE(X|A\wedge B)\in L(X|A\wedge B)=]-\infty,\LE(X|A\wedge B)[$.

Finally, we prove that $\LE(X|B)=\LP(X|B)\ \forall X|B\in\dset$.
If $X|B\in\dset$, taking $n=1$, $s_{1}=1$, $X_{1}|B_{1}=X|B$ in the definition of $\LE(X|B)$,
$\sup \{B(X-\LP(X|B))-B(X-\alpha)|B\}=\alpha-\LP(X|B)<0$, $\forall \alpha<\LP(X|B)$.
Hence $\LE(X|B)\geq\LP(X|B)$.
For the reverse inequality, note  that $\forall X|B\in\dset$, $\forall X_{i}|B_{i}\in\dset$, $\forall s_{i}\geq 0\ (i=1,\ldots,n)$,
$\sup\{\sum_{i=1}^{n}g_{i}-B(X-\LP(X|B))|\bigvee_{i=1}^{n}B_{i}\vee B\}\geq 0$,
by the coherence of $\LP$ on $\dset$. It ensues $\LP(X|B)\notin L(X|B)=]-\infty,\LE(X|B)[$.
\fineproof

Theorem \ref{thm:extension_conditional_lower_prevision} lets us extend any W-coherent conditional lower prevision
$\LP:\dset\rightarrow\rset$ to any set $\dset^{\prime}$($\supset\dset$) which meets the structure requirements of $\dset^{*}$
in Theorem \ref{thm:conditional coherent prevision on linear space}.
The set $\dset^{\prime}$ does not necessarily satisfy these requirements.
When it does not,
consider a partition $\bset$ on which the random variables in $\dset^{\prime}$ are defined
and let $\xset$ be the set of all random variables on $\bset$, $\eset^{\emptyset}$ and $\dset^{*}$
as in Theorem \ref{thm:conditional coherent prevision on linear space}.
By Theorem \ref{thm:extension_conditional_lower_prevision},
$\LE:\dset^{*}\rightarrow\rset$ is a W-coherent conditional lower previsions,
extending $\LP$ to $\dset^{*}$ and therefore to $\dset^{\prime}\subset\dset^{*}$ as well.


\section*{Acknowledgements}
We are grateful to the referees for their constructive suggestions.

*NOTICE: This is the authors' version of a work that was accepted for publication in the International Journal of Approximate Reasoning. Changes resulting from the publishing process, such as peer review, editing, corrections, structural formatting, and other quality control mechanisms may not be reflected in this document. Changes may have been made to this work since it was submitted for publication. A definitive version was subsequently published in the International Journal of Approximate Reasoning, vol. 50, issue 4, April 2009, doi:10.1016/j.ijar.2008.10.002

$\copyright$ Copyright Elsevier

http://www.sciencedirect.com/science/article/pii/S0888613X08001643

\end{document}